\documentclass[11pt, a4paper]{article}
\usepackage[cp1251]{inputenc}
\usepackage[russian]{babel}
\usepackage{amsmath} \usepackage{euscript}
\usepackage{mymatrx5}
\def\mymatrix{\MyMatrixwithdelims..}
\Linestrue\Autonumtrue
\usepackage{mathrsfs}
\usepackage{amssymb}
\usepackage[matrix, arrow, curve]{xy}
\begin{document}

\setcounter{MaxMatrixCols}{14}
\newcounter{bnomer} \newcounter{snomer}
\newcounter{bsnomer}
\setcounter{bnomer}{0}
\renewcommand{\thesnomer}{\thebnomer.\arabic{snomer}}
\renewcommand{\thebsnomer}{\thebnomer.\arabic{bsnomer}}
\renewcommand{\refname}{\begin{center}\large{\textbf{References}}\end{center}}

\newcommand{\sect}[1]{%
\setcounter{snomer}{0}\setcounter{bsnomer}{0}
\refstepcounter{bnomer}
\par\bigskip\begin{center}\large{\textbf{\arabic{bnomer}. {#1}}}\end{center}}
\newcommand{\sst}[1]{%
\refstepcounter{bsnomer}
\par\bigskip\textbf{\arabic{bnomer}.\arabic{bsnomer}. {#1}.}}
\newcommand{\defi}[1]{%
\refstepcounter{snomer}
\par\medskip\textbf{Definition \arabic{bnomer}.\arabic{snomer}. }{#1}\par\medskip}
\newcommand{\theo}[2]{%
\refstepcounter{snomer}
\par\textbf{Теорема \arabic{bnomer}.\arabic{snomer}. }{#2} {\emph{#1}}\hspace{\fill}$\square$\par}
\newcommand{\mtheop}[2]{%
\refstepcounter{snomer}
\par\textbf{Theorem \arabic{bnomer}.\arabic{snomer}. }{\emph{#1}}
\par\textsc{Proof}. {#2}\hspace{\fill}$\square$\par}
\newcommand{\mcorop}[2]{%
\refstepcounter{snomer}
\par\textbf{Corollary \arabic{bnomer}.\arabic{snomer}. }{\emph{#1}}
\par\textsc{Proof}. {#2}\hspace{\fill}$\square$\par}
\newcommand{\mtheo}[1]{%
\refstepcounter{snomer}
\par\medskip\textbf{Theorem \arabic{bnomer}.\arabic{snomer}. }{\emph{#1}}\par\medskip}
\newcommand{\mlemm}[1]{%
\refstepcounter{snomer}
\par\medskip\textbf{Lemma \arabic{bnomer}.\arabic{snomer}. }{\emph{#1}}\par\medskip}
\newcommand{\mprop}[1]{%
\refstepcounter{snomer}
\par\medskip\textbf{Proposition \arabic{bnomer}.\arabic{snomer}. }{\emph{#1}}\par\medskip}
\newcommand{\theobp}[2]{%
\refstepcounter{snomer}
\par\textbf{Теорема \arabic{bnomer}.\arabic{snomer}. }{#2} {\emph{#1}}\par}
\newcommand{\theop}[2]{%
\refstepcounter{snomer}
\par\textbf{Theorem \arabic{bnomer}.\arabic{snomer}. }{\emph{#1}}
\par\textsc{Proof}. #2\hspace{\fill}$\square$\par}
\newcommand{\theosp}[2]{%
\refstepcounter{snomer}
\par\textbf{Теорема \arabic{bnomer}.\arabic{snomer}. }{\emph{#1}}
\par\textbf{Схема доказательства}. {#2}\hspace{\fill}$\square$\par}
\newcommand{\exam}[1]{%
\refstepcounter{snomer}
\par\medskip\textbf{Example \arabic{bnomer}.\arabic{snomer}. }{#1}\par\medskip}
\newcommand{\deno}[1]{%
\refstepcounter{snomer}
\par\textbf{Definition \arabic{bnomer}.\arabic{snomer}. }{#1}\par}
\newcommand{\post}[1]{%
\refstepcounter{snomer}
\par\textbf{Предложение \arabic{bnomer}.\arabic{snomer}. }{\emph{#1}}\hspace{\fill}$\square$\par}
\newcommand{\postp}[2]{%
\refstepcounter{snomer}
\par\medskip\textbf{Proposition \arabic{bnomer}.\arabic{snomer}. }{\emph{#1}}%
\ifhmode\par\fi\textsc{Proof}. {#2}\hspace{\fill}$\square$\par\medskip}
\newcommand{\lemm}[1]{%
\refstepcounter{snomer}
\par\textbf{Lemma \arabic{bnomer}.\arabic{snomer}. }{\emph{#1}}\hspace{\fill}$\square$\par}
\newcommand{\lemmp}[2]{%
\refstepcounter{snomer}
\par\medskip\textbf{Lemma \arabic{bnomer}.\arabic{snomer}. }{\emph{#1}}
\par\textsc{Proof}. #2\hspace{\fill}$\square$\par\medskip}
\newcommand{\coro}[1]{%
\refstepcounter{snomer}
\par\textbf{Corollary \arabic{bnomer}.\arabic{snomer}. }{\emph{#1}}\hspace{\fill}$\square$\par}
\newcommand{\mcoro}[1]{%
\refstepcounter{snomer}
\par\textbf{Corollary \arabic{bnomer}.\arabic{snomer}. }{\emph{#1}}\par\medskip}
\newcommand{\corop}[2]{%
\refstepcounter{snomer}
\par\textbf{Corollary \arabic{bnomer}.\arabic{snomer}. }{\emph{#1}}
\par\textsc{Proof}. {#2}\hspace{\fill}$\square$\par}
\newcommand{\nota}[1]{%
\refstepcounter{snomer}
\par\medskip\textbf{Remark \arabic{bnomer}.\arabic{snomer}. }{#1}\par\medskip}
\newcommand{\propp}[2]{%
\refstepcounter{snomer}
\par\medskip\textbf{Proposition \arabic{bnomer}.\arabic{snomer}. }{\emph{#1}}
\par\textsc{Proof}. {#2}\hspace{\fill}$\square$\par\medskip}
\newcommand{\hypo}[1]{%
\refstepcounter{snomer}
\par\medskip\textbf{Conjecture \arabic{bnomer}.\arabic{snomer}. }{\emph{#1}}\par\medskip}
\newcommand{\prop}[1]{%
\refstepcounter{snomer}
\par\textbf{Proposition \arabic{bnomer}.\arabic{snomer}. }{\emph{#1}}\hspace{\fill}$\square$\par}

\newcommand{\Ind}[3]{%
\mathrm{Ind}_{#1}^{#2}{#3}}
\newcommand{\Res}[3]{%
\mathrm{Res}_{#1}^{#2}{#3}}
\newcommand{\epsi}{\epsilon}
\newcommand{\tri}{\triangleleft}
\newcommand{\Supp}[1]{%
\mathrm{Supp}(#1)}

\makeatletter
\def\iddots{\mathinner{\mkern1mu\raise\p@
\vbox{\kern7\p@\hbox{.}}\mkern2mu
\raise4\p@\hbox{.}\mkern2mu\raise7\p@\hbox{.}\mkern1mu}}
\makeatother

\newcommand{\reg}{\mathrm{reg}}
\newcommand{\empr}[2]{[-{#1},{#1}]\times[-{#2},{#2}]}
\newcommand{\sreg}{\mathrm{sreg}}
\newcommand{\codim}{\mathrm{codim}\,}
\newcommand{\chara}{\mathrm{char}\,}
\newcommand{\rk}{\mathrm{rk}\,}
\newcommand{\chr}{\mathrm{ch}\,}
\newcommand{\id}{\mathrm{id}}
\newcommand{\Ad}{\mathrm{Ad}}
\newcommand{\Ker}{\mathrm{Ker}\,}
\newcommand{\End}{\mathrm{End}}
\newcommand{\Ann}{\mathrm{Ann}\,}
\newcommand{\col}{\mathrm{col}}
\newcommand{\row}{\mathrm{row}}
\newcommand{\low}{\mathrm{low}}
\newcommand{\sur}{\twoheadrightarrow}
\newcommand{\inj}{\hookrightarrow}
\newcommand{\pho}{\hphantom{\quad}\vphantom{\mid}}
\newcommand{\fho}[1]{\vphantom{\mid}\setbox0\hbox{00}\hbox to \wd0{\hss\ensuremath{#1}\hss}}
\newcommand{\wt}{\widetilde}
\newcommand{\wh}{\widehat}
\newcommand{\ad}[1]{\mathrm{ad}_{#1}}
\newcommand{\tr}{\mathrm{tr}\,}
\newcommand{\GL}{\mathrm{GL}}
\newcommand{\Pf}{\mathrm{Pf}}
\newcommand{\Prim}{\mathrm{Prim}\,}
\newcommand{\Cent}{\mathrm{Cent}\,}
\newcommand{\SL}{\mathrm{SL}}
\newcommand{\SO}{\mathrm{SO}}
\newcommand{\Sp}{\mathrm{Sp}}
\newcommand{\Mat}{\mathrm{Mat}}
\newcommand{\Stab}{\mathrm{Stab}}
\newcommand{\ilm}{\varinjlim}
\newcommand{\lee}{\leq}
\newcommand{\gee}{\geq}

\newcommand{\vfi}{\varphi}
\newcommand{\vpi}{\varpi}
\newcommand{\teta}{\vartheta}
\newcommand{\Bfi}{\Phi}
\newcommand{\Fp}{\mathbb{F}}
\newcommand{\Rp}{\mathbb{R}}
\newcommand{\Zp}{\mathbb{Z}}
\newcommand{\Cp}{\mathbb{C}}
\newcommand{\Np}{\mathbb{N}}
\newcommand{\Xt}{\mathfrak{X}}
\newcommand{\ut}{\mathfrak{u}}
\newcommand{\at}{\mathfrak{a}}
\newcommand{\nt}{\mathfrak{n}}
\newcommand{\mt}{\mathfrak{m}}
\newcommand{\htt}{\mathfrak{h}}
\newcommand{\spt}{\mathfrak{sp}}
\newcommand{\slt}{\mathfrak{sl}}
\newcommand{\ot}{\mathfrak{so}}
\newcommand{\rt}{\mathfrak{r}}
\newcommand{\rad}{\mathfrak{rad}}
\newcommand{\bt}{\mathfrak{b}}
\newcommand{\gt}{\mathfrak{g}}
\newcommand{\vt}{\mathfrak{v}}
\newcommand{\pt}{\mathfrak{p}}
\newcommand{\kt}{\mathfrak{k}}
\newcommand{\Po}{\mathcal{P}}
\newcommand{\Uo}{\EuScript{U}}
\newcommand{\Fo}{\EuScript{F}}
\newcommand{\Do}{\EuScript{D}}
\newcommand{\Eo}{\EuScript{E}}
\newcommand{\Iu}{\mathcal{I}}
\newcommand{\Du}{\mathcal{D}}
\newcommand{\Uu}{\mathcal{U}}
\newcommand{\Mo}{\mathcal{M}}
\newcommand{\Nu}{\mathcal{N}}
\newcommand{\Ro}{\mathcal{R}}
\newcommand{\Co}{\mathcal{C}}
\newcommand{\Lo}{\mathcal{L}}
\newcommand{\Ou}{\mathcal{O}}
\newcommand{\Au}{\mathcal{A}}
\newcommand{\Vu}{\mathcal{V}}
\newcommand{\Bu}{\mathcal{B}}
\newcommand{\Sy}{\mathcal{Z}}
\newcommand{\Sb}{\mathcal{F}}
\newcommand{\Gr}{\mathcal{G}}
\newcommand{\rtc}[1]{C_{#1}^{\mathrm{red}}}

\author{Mikhail V. Ignatyev}

\date{}
\title{\Large{Centrally generated primitive ideals of $U(\nt)$ in types $B$ and $D$\mbox{$\vphantom{1}$}\footnotetext{The work on Section~\ref{sect:infinite_dim} was performed at the NRU HSE with the support from the Russian Science Foundation, grant no. 16--41--01013. The work on Section~\ref{sect:finite_dim} has been supported by RFBR grant no. 16--01--00154a and by ISF grant no. 797/14.}}
}\maketitle

\begin{center}
\begin{tabular}{p{15cm}}
\small{\textsc{Abstract}. We study the centrally generated primitive ideals of $U(\nt)$, where $\nt$ is the (locally) nilpotent radical of a (splitting) Borel subalgebra of a simple complex Lie algebra $\gt=\ot_{2n+1}(\Cp)$, $\ot_{2n}(\Cp)$, $\ot_{\infty}(\Cp)$. In the infinite-dimensional setting, there are infinitely many isomorphism classes of Lie algebras $\nt$, and we fix $\nt$ with ``largest possible'' center of $U(\nt)$. We characterize the centrally generated primitive ideals of $U(\nt)$ in terms of the Dixmier map and the Kostant cascade.}\\\\
\small{\textbf{Keywords:} Dixmier map, Kostant cascade, center of enveloping algebra, centrally generated primitive ideal, locally nilpotent Lie algebra.}\\
\small{\textbf{AMS subject classification:} 17B65, 17B35, 17B10, 17B08.}
\end{tabular}
\end{center}

\sect{Introduction}

The theory of primitive ideals in enveloping algebras of Lie algebras has its roots in the re\-pre\-sen\-tation theory of Lie algebras. However, classifying irreducible representations of Lie algebras is not feasible except in few very special cases, while a classification of annihilators of irreducible representations, i.e., of primitive ideals, can be achieved in much greater generality. This idea goes back to J. Dixmier and his seminar, and for semisimple or solvable finite-dimensional Lie algebras there is an extensive theory of primitive ideals.

In the case when $\nt$ is a finite-dimensional nilpotent Lie algebra, the primitive ideals in the universal enveloping algebra $U(\nt)$ can be described in terms of the Dixmier map assigning to any linear form $f\in\nt^*$ a primitive ideal $J(f)$ of $U(\nt)$. If $\nt$ is abelian, $J(f)$ is simply the annihilator of $f$. For a general finite-dimensional nilpotent Lie algebra $\nt$, the theory of primitive ideals retains many properties from the abelian case: in particular, $J(f)$ is always a maximal ideal and every primitive ideal in $U(\nt)$ is of the form $J(f)$ for some $f\in\nt^*$. Moreover, $J(f)=J(f')$ if and only if $f$ and $f'$ belong to the same coadjoint orbit in $\nt^*$. Note that, according to the Kostant--Kirillov orbit method, the set of all unitary irreducible representations of the Lie group $\exp\nt$ is also in one-to-one correspondence with the set of all coadjoint orbits in $\nt^*$.

Suppose that $\nt$ is the nilradical of a Borel subalgebra $\bt$ of a complex finite-dimensional semisimple Lie algebra $\gt$. The description of the center $Z(\nt)$ of $U(\nt)$ goes back to J. Dixmier, A. Joseph and B.~Kostant. It turns out that $Z(\nt)$ is a polynomial algebra whose generators are parametrized by the positive roots from the Kostant cascade $\Bu$, a certain strongly orthogonal subset of the set $\Phi^+$ of positive roots of $\gt$ with respect to $\bt$ (see Subsection~\ref{sst:center_fd} for the precise definition). Let $\{e_{\alpha},~\alpha\in\Phi^+\}$ be a basis of $\nt$ consisting of root vectors. We say that a linear form $f\in\nt^*$ is a Kostant form if $f(e_{\alpha})=0$ for $\alpha\notin\Bu$ and $f(e_{\beta})\neq0$ for $\beta\in\Bu\setminus\Delta$, where $\Delta$ is the set of simple roots. Note that the coadjoint orbit of a Kostant form has maximal possible dimension. By definition, a primitive ideal $J$ of $U(\nt)$ is centrally generated if it is generated as an ideal by its intersection with $Z(\nt)$. It was proved in \cite[Theorem~3.1]{IgnatyevPenkov1} that, for $\Phi=A_{n-1}$ or $\Phi=C_n$, $J$ is centrally generated if and only if $J=J(f)$ for a certain Kostant form~$f$. In this paper, we prove that this fact is also true for $\Phi=B_n$ or $\Phi=D_n$.

More precisely, let $P_{\beta}$, $\beta\in\Bu$, be the set of canonical generators of $Z(\nt)$ (the explicit formulas for~$P_{\beta}$ are given in Subsection~\ref{sst:center_fd}). Let $J$ be a primitive ideal of $U(\nt)$. Since $J$ is the annihilator of a~simple $\nt$-module, given $\beta\in\Bu$, there exists the unique $c_{\beta}\in\Cp$ such that $P_{\beta}-c_{\beta}\in J$. Our first main result (see Theorem~\ref{theo:ideal_fd_Bn_Dn}) claims that the following conditions are equivalent:
\begin{equation*}
\begin{split}
&\text{\textup{i)} $J$ is centrally generated};\\
&\text{\textup{ii)} all scalars $c_{\beta}$, $\beta\in\Bu\setminus\Delta$, are nonzero};\\
&\text{\textup{iii)} $J=J(f)$ for a Kostant form $f$}.
\end{split}
\end{equation*}
If these conditions are satisfied, then we present an explicit way how to reconstruct $f$ by $J$.

Suppose now that $\nt$ is the locally nilpotent radical of a Borel subalgebra $\bt$ of the simple finitary infinite-dimensional complex Lie algebra $\gt=\ot_{\infty}(\Cp)$. The center of the enveloping algebra $U(\nt)$ was described in \cite[Theorem~2.6]{IgnatyevPenkov1}. There are infinitely many isomorphism classes of Lie algebras $\nt$, and we concentrate on the case when $U(\nt)$ has ``largest possible'' center, which in technical terms means that $\nt$ admits a filtration by finite-dimensional subalgebras so that the centers of their enveloping algebras are appropriately aligned. This latter requirement singles out only one isomorphism class of subalgebras~$\nt$. For such $\nt$ we construct a partial Dixmier map defined for certain linear forms
$f\in\nt^*$, closely related to the Kostant cascade of $\nt$. We refer to these forms as Kostant forms. Our second main result (Theorem~\ref{theo:ideal_ifd}) implies then that, similarly to the finite-dimensional case, the Dixmier map establishes a one-to-one correspondence between Kostant forms and centrally generated primitive ideals in $U(\nt)$. This provides an explicit description of the centrally generated primitive ideals of $U(\nt)$.

The paper is organized as follows. Section~\ref{sect:finite_dim} is devoted to the finite-dimensional case. In Sub\-sec\-tion~\ref{sst:center_fd}, we briefly recall the Kostant's characterization of $Z(\nt)$ and present an explicit description of the canonical generators of $Z(\nt)$ from \cite[Subsection 2.1]{IgnatyevPenkov1} based on A. Panov's work \cite{Panov2}. (For the reader's convenience, we present the answer for all classical Lie algebras $\gt$.) Subsection~\ref{sst:ideals_fd} contains the proof of our first main result, Theorem~\ref{theo:ideal_fd_Bn_Dn}. The key step in the proof is Proposition~\ref{prop:J_c_primitive_Bn_Dn}. We also find explicitly a set of ``canonical'' generators of the quotient algebra $U(\nt)/J$, which is isomorphic to a Weyl algebra, see Theorem~\ref{theo:Weyl_quotient}. Next, in Section~\ref{sect:infinite_dim}, we consider the infinite-dimensional case. In Subsection~\ref{sst:centers_ifd}, we briefly recall a description of $Z(\nt)$ from \cite[Theorem 2.6]{IgnatyevPenkov1}. Finally, in Subsection~\ref{sst:ideals_ifd} we restrict ourselves on the special choice of $\nt$ and, using Proposition~\ref{prop:n_annihilating_n_plus_1}, prove our second main result, Theorem~\ref{theo:ideal_ifd}, describing centrally generated ideals of $U(\nt)$.

\textsc{Acknowledgments}. I express my gratitude to A. Melnikov, A. Panov, I. Penkov and A. Shev\-chen\-ko for useful discussions. The work on Section~\ref{sect:finite_dim} was done during my stay at the University of Haifa. I thank this institution for its hospitality.

\sect{Finite-dimensional\label{sect:finite_dim} case}

\sst{The\label{sst:center_fd} center of $U(\nt)$} Let $n\in\Zp_{>0}$. Throughout this subsection $\gt$ denotes one of the Lie algebras $\slt_n(\Cp)$, $\ot_{2n}(\Cp)$, $\ot_{2n+1}(\Cp)$ or $\spt_{2n}(\Cp)$. The algebra $\ot_{2n}(\Cp)$ (respectively, $\ot_{2n+1}(\Cp)$ and $\spt_{2n}(\Cp)$) is realized as the sub\-al\-gebra of $\slt_{2n}(\Cp)$ (respectively, $\slt_{2n+1}(\Cp)$ and $\slt_{2n}(\Cp)$) consisting of all~$x$ such that $\beta(u,xv)+\beta(xu,v)=0$ for all $u,v$ in $\Cp^{2n}$ (respectively, in $\Cp^{2n+1}$ and $\Cp^{2n}$), where
\begin{equation*}
\beta(u,v)=\begin{cases}
\sum\nolimits_{i=1}^n(u_iv_{-i}+u_{-i}v_i)&\text{for }\ot_{2n}(\Cp),\\
u_0v_0+\sum\nolimits_{i=1}^n(u_iv_{-i}+u_{-i}v_i)&\text{for }\ot_{2n+1}(\Cp),\\
\sum\nolimits_{i=1}^n(u_iv_{-i}-u_{-i}v_i)&\text{for }\spt_{2n}(\Cp).
\end{cases}
\end{equation*}
Here for $\ot_{2n}(\Cp)$ (respectively, for $\ot_{2n+1}$ and $\spt_{2n}(\Cp))$ we denote by $e_1,\ldots,e_n,e_{-n},\ldots,e_{-1}$ (res\-pec\-tively, by $e_1,\ldots,e_n,e_0,e_{-n},\ldots,e_{-1}$ and $e_1,\ldots,e_n,e_{-n},\ldots,e_{-1}$) the standard basis of $\Cp^{2n}$ (res\-pec\-tively, of $\Cp^{2n+1}$ and $\Cp^{2n}$), and by $x_i$ the coordinate of a vector $x$ corresponding to $e_i$.

The set of all diagonal matrices from $\gt$ is a Cartan subalgebra of $\gt$; we denote it by $\htt$. Let $\Phi$ be the root system of $\gt$ with respect to $\htt$. Note that $\Phi$ is of type $A_{n-1}$ (respectively, $D_n$, $B_n$ and $C_n$) for $\slt_n(\Cp)$ (respectively, for $\ot_{2n}(\Cp)$, $\ot_{2n+1}(\Cp)$ and $\spt_{2n}(\Cp)$). The set of all upper-triangular matrices from $\gt$ is a Borel subalgebra of $\gt$ containing $\htt$; we denote it by $\bt$. Let $\Phi^+$ be the set of positive roots with respect to $\bt$. As usual, we identify $\Phi^+$ with the following subset of $\Rp^n$:
\begin{equation*}\predisplaypenalty=0
\begin{split}
A_{n-1}^+&=\{\epsi_i-\epsi_j,~1\leq i<j\leq n\},\\
B_n^+&=\{\epsi_i-\epsi_j,~1\leq i<j\leq n\}\cup\{\epsi_i+\epsi_j,~1\leq i<j\leq n\}\cup\{\epsi_i,~1\leq i\leq n\},\\
C_n^+&=\{\epsi_i-\epsi_j,~1\leq i<j\leq n\}\cup\{\epsi_i+\epsi_j,~1\leq i<j\leq n\}\cup\{2\epsi_i,~1\leq i\leq n\},\\
D_n^+&=\{\epsi_i-\epsi_j,~1\leq i<j\leq n\}\cup\{\epsi_i+\epsi_j,~1\leq i<j\leq n\}.\\
\end{split}
\end{equation*}
Here $\{\epsi_i\}_{i=1}^n$ is the standard basis of $\Rp^n$.

Denote by $\nt$ the algebra of all strictly upper-triangular matrices from $\gt$. Then $\nt$ has a basis consisting of root vectors $e_{\alpha}$, $\alpha\in\Phi^+$, where
\begin{equation*}\predisplaypenalty=0
\begin{split}
e_{\epsi_i}&=\sqrt{2}(e_{0,i}-e_{-i,0}),~e_{2\epsi_i}=e_{i,-i},\\
e_{\epsi_i-\epsi_j}&=\begin{cases}
e_{i,j}&\text{for }A_{n-1},\\
e_{i,j}-e_{-j,-i}&\text{for }B_n,~C_n\text{ and }D_n,
\end{cases}\\
e_{\epsi_i+\epsi_j}&=\begin{cases}
e_{i,-j}-e_{j,-i}&\text{for }B_n\text{ and }D_n,\\
e_{i,-j}+e_{j,-i}&\text{for }C_n,
\end{cases}
\end{split}
\end{equation*}
and $e_{i,j}$ are the usual elementary matrices. For $\ot_{2n}(\Cp)$ (respectively, for $\ot_{2n+1}(\Cp)$ and $\spt_{2n}(\Cp)$) we index the rows (from left to right) and the columns (from top to bottom) of matrices by the numbers $1,\ldots,n,-n,\ldots,-1$ (respectively, by the numbers $1,\ldots,n,0,-n,\ldots,-1$ and $1,\ldots,n,-n,\ldots,-1$). Note that $\gt=\htt\oplus\nt\oplus\nt_-$, where $\nt_-=\langle e_{-\alpha},~\alpha\in\Phi^+\rangle_{\Cp}$, and, by definition, $e_{-\alpha}=e_{\alpha}^T$. (The superscript~$T$ always indicates matrix transposition.) The set $\{e_{\alpha},~\alpha\in\Phi\}$ can be extended to a unique Chevalley basis of $\gt$.

Let $G$ be one of the following classical Lie groups: $\SL_n(\Cp)$, $\SO_{2n}(\Cp)$, $\SO_{2n+1}(\Cp)$ or $\Sp_{2n}(\Cp)$. The group $\SO_{2n}(\Cp)$ (respectively, $\SO_{2n+1}$ and $\Sp_{2n}(\Cp)$) is realized as the subgroup of $\SL_{2n}(\Cp)$ (respectively, of $\SL_{2n+1}(\Cp)$ and $\SL_{2n}(\Cp)$) which preserves the form $\beta$. Let $H$ (respectively, $B$ and $N$) be the set of all diagonal (respectively, upper-triangular and upper-triangular with 1 on the diagonal) matrices from~$G$. Then $H$ is a maximal torus of $G$, $B$ is a Borel subgroup of $G$ containing $H$, $N$ is the unipotent radical of $B$, and $\gt$ (respectively, $\htt$, $\bt$ and $\nt$) is the Lie algebra of $G$ (respectively, of $H$, $B$ and $N$).

Denote by $U(\nt)$ the enveloping algebra of $\nt$, and by $S(\nt)$ the symmetric algebra of $\nt$. Then $\nt$ and $S(\nt)$ are $B$-modules as $B$ normalizes $N$. Denote by $Z(\nt)$ the center of $U(\nt)$. It is well-known that the restriction of the symmetrization map $$\sigma\colon S(\nt)\to U(\nt),~x^k\mapsto x^k,~x\in\nt,~k\in\Zp_{\geq0},$$ to the algebra $S(\nt)^N$ of $N$-invariants is an algebra isomorphism between $S(\nt)^N$ and $Z(\nt)$.

Denote by $\Bu$ the following subset of $\Phi^+$:
\begin{equation*}
\Bu=\begin{cases}\bigcup\nolimits_{1\leq i\leq[n/2]}\{\epsi_i-\epsi_{n-i+1}\}&\text{for }A_{n-1},\\
\bigcup\nolimits_{1\leq i\leq n/2}\{\epsi_{2i-1}-\epsi_{2i},~\epsi_{2i-1}+\epsi_{2i}\}&\text{for }B_n,~n\text{ even},\\
\bigcup\nolimits_{1\leq i\leq[n/2]}\{\epsi_{2i-1}-\epsi_{2i},~\epsi_{2i-1}+\epsi_{2i}\}\cup\{\epsi_n\}&\text{for }B_n,~n\text{ odd},\\
\bigcup\nolimits_{1\leq i\leq n}\{2\epsi_i\}&\text{for }C_n,\\
\bigcup\nolimits_{1\leq i\leq[n/2]}\{\epsi_{2i-1}-\epsi_{2i},~\epsi_{2i-1}+\epsi_{2i}\}&\text{for }D_n.\\
\end{cases}
\end{equation*}
Note that $\Bu$ is a maximal strongly orthogonal subset of $\Phi^+$, i.e., $\Bu$ is maximal with the property that if $\alpha,\beta\in\Bu$ then neither $\alpha-\beta$ nor $\alpha+\beta$ belongs to $\Phi^+$. We call $\Bu$ the \emph{Kostant cascade} of orthogonal roots in~$\Phi^+$.

We next present a canonical set of generators of $Z(\nt)$ (or, equivalently, of $S(\nt)^N$), whose description goes back to J. Dixmier, A. Joseph and B. Kostant \cite{Dixmier1}, \cite{Joseph1}, \cite{Kostant1}, \cite{Kostant2}. We can consider $\Zp\Phi$, the $\Zp$-linear span of $\Phi$, as a subgroup of the group $\Xt$ of rational multiplicative characters of $H$ by putting $\pm\epsi_i(h)=h_{i,i}^{\pm1}$, where $h_{i,i}$ is the $i$-th diagonal element of a matrix $h\in H$. Recall that a vector $\lambda\in\Rp^n$ is called a \emph{weight} of $H$ if $c(\alpha,\lambda)=2(\alpha,\lambda)/(\alpha,\alpha)$ is an integer for any $\alpha\in\Phi^+$, where $(\cdot,\cdot)$ is the standard inner product on $\Rp^n$. A weight $\lambda$ is called \emph{dominant} if $c(\alpha,\lambda)\geq0$ for all $\alpha\in\Phi^+$. An element $a$ of an $H$-module is called an \emph{$H$-weight vector}, if there exists $\nu\in\Xt$ such that $h\cdot a=\nu(h)a$ for all $h\in H$. By \cite[Theorems 6, 7]{Kostant2}, every $H$-weight occurs in $S(\nt)^N$ with multiplicity at most 1. Furthermore, there exist unique (up to scalars) prime polynomials $\xi_{\beta}\in S(\nt)^N$, $\beta\in\Bu$, such that each $\xi_{\beta}$ is an $H$-weight polynomial of a dominant weight $\vpi_{\beta}$ belonging to the $\Zp$-linear span $\Zp\Bu$ of $\Bu$. A remarkable fact is that
\begin{equation}
\xi_{\beta},~\beta\in\Bu,\text{ are algebraically independent generators of }S(\nt)^N,\label{formula:xi_i}
\end{equation}
so $S(\nt)^N$ and $Z(\nt)$ are polynomial rings. We call $\{\xi_{\beta},~\beta\in\Bu\}$ the set of \emph{canonical generators} of $S(\nt)^N$. Put $m=|\Bu|$. It turns out that the weights $\vpi_{\beta}$'s have the following form \cite[Theorem 2.12]{Panov2}.
\begin{equation}\label{formula:weights_of_xi_i}\text{\begin{tabular}{|l|l|}
\hline
$\Phi=A_{n-1}$&$\vpi_{\beta}=2\epsi_1+\ldots+2\epsi_{i-1}+\epsi_i$ for $\beta=\epsi_i-\epsi_{n-i+1}$, $1\leq i\leq m$\\
\hline
$\Phi=B_n$&$\vpi_{\beta}=\begin{cases}2\epsi_1+\ldots+2\epsi_i&\text{for $\beta=\epsi_i-\epsi_{i+1}$ with odd $i<m$},\\
\epsi_1+\ldots+\epsi_i&\text{otherwise}
\end{cases}$\\
\hline
$\Phi=C_n$&$\vpi_{\beta}=2\epsi_1+\ldots+2\epsi_i$ for $\beta=2\epsi_i$, $1\leq i\leq m$\\
\hline
$\Phi=D_n$&$\vpi_{\beta}=\begin{cases}2\epsi_1+\ldots+2\epsi_i&\text{for $\beta=\epsi_i-\epsi_{i+1}$ with odd $i<m-1$},\\
2\epsi_1+\ldots+2\epsi_{n-1}&\text{for $\beta=\epsi_{n-2}-\epsi_{n-1}$ when $n$ is odd},\\
\epsi_1+\ldots+\epsi_{n-1}-\epsi_n&\text{for $\beta=\epsi_{n-1}-\epsi_n$ when $n$ is even},\\
\epsi_1+\ldots+\epsi_i&\text{otherwise}
\end{cases}$\\
\hline
\end{tabular}}
\end{equation}

The following description of $\xi_{\beta}$ for classical root systems, given in \cite[Subsection 2.1]{IgnatyevPenkov1}, follows from \cite{Panov1} (see also \cite{LipsmanWolf1} and \cite{Joseph-Faucaunt-Millet1}). Our notation here slightly differs from the one used in \cite{IgnatyevPenkov1}.

i) $\Phi=A_{n-1}$. Here, for $1\leq i\leq m=[n/2]$,
\begin{equation}
\xi_{\epsi_i-\epsi_{n-i+1}}=\begin{vmatrix}\label{formula:Delta_i_A_n}
e_{1,n-i+1}&\ldots&e_{1,n-1}&e_{1,n}\\
e_{2,n-i+1}&\ldots&e_{2,n-1}&e_{2,n}\\
\vdots&\iddots&\vdots&\vdots\\
e_{i,n-i+1}&\ldots&e_{i,n-1}&e_{i,n}\\
\end{vmatrix}=\begin{vmatrix}
e_{\epsi_1-\epsi_{n-i+1}}&\ldots&e_{\epsi_1-\epsi_{n-1}}&e_{\epsi_1-\epsi_n}\\
e_{\epsi_2-\epsi_{n-i+1}}&\ldots&e_{\epsi_2-\epsi_{n-1}}&e_{\epsi_2-\epsi_n}\\
\vdots&\iddots&\vdots&\vdots\\
e_{\epsi_i-\epsi_{n-i+1}}&\ldots&e_{\epsi_i-\epsi_{n-1}}&e_{\epsi_i-\epsi_n}\\
\end{vmatrix}.
\end{equation}

ii) $\Phi=C_n$. Here, for $1\leq i\leq m=n$,
\begin{equation}
\xi_{2\epsi_i}=\begin{vmatrix}
e_{\epsi_1+\epsi_i}&\ldots&e_{\epsi_1+\epsi_3}&e_{\epsi_1+\epsi_2}&2e_{2\epsi_1}\\
e_{\epsi_2+\epsi_i}&\ldots&e_{\epsi_2+\epsi_3}&2e_{2\epsi_2}&e_{\epsi_1+\epsi_2}\\
e_{\epsi_3+\epsi_i}&\ldots&2e_{2\epsi_3}&e_{\epsi_2+\epsi_3}&e_{\epsi_1+\epsi_3}\\
\vdots&\iddots&\vdots&\vdots&\vdots\\
2e_{2\epsi_i}&\ldots&e_{\epsi_3+\epsi_i}&e_{\epsi_2+\epsi_i}&e_{\epsi_1+\epsi_i}\\
\end{vmatrix}.\label{formula:Delta_i_C_n}
\end{equation}

iii) $\Phi=D_n$. If $i$ is odd, then
\begin{equation}
\xi_{\epsi_i+\epsi_{i+1}}^2=\pm\begin{vmatrix}
e_{\epsi_1+\epsi_{i+1}}&\ldots&e_{\epsi_1+\epsi_3}&e_{\epsi_1+\epsi_2}&0\\
e_{\epsi_2+\epsi_{i+1}}&\ldots&e_{\epsi_2+\epsi_3}&0&-e_{\epsi_1+\epsi_2}\\
e_{\epsi_3+\epsi_{i+1}}&\ldots&0&-e_{\epsi_2+\epsi_3}&-e_{\epsi_1+\epsi_3}\\
\vdots&\iddots&\vdots&\vdots&\vdots\\
0&\ldots&-e_{\epsi_3+\epsi_{i+1}}&-e_{\epsi_2+\epsi_{i+1}}&-e_{\epsi_1+\epsi_{i+1}}\\
\end{vmatrix}.\label{formula:Delta_i_D_n_pf}\postdisplaypenalty=10000
\end{equation}
(After suitable reordering of indices, the matrix in the right-hand side becomes skew-symmetric, so $\xi_{\epsi_i+\epsi_{i+1}}$ is nothing but its Pfaffian.) Our normalization is such that the term $e_{\epsi_1+\epsi_2}e_{\epsi_3+\epsi_4}\ldots e_{\epsi_i+\epsi_{i+1}}$ enters $\xi_{\epsi_i+\epsi_{i+1}}$ with coefficient 1.

Next, assume that $i\leq n-1$ is even. Let $\Uu^i$ be the $(i\times i)$-matrix with entries from $S(\nt)$ defined by $(\Uu^i)_{a,b}=-(\Uu^i)_{i-b+1,i-a+1}=e_{\epsi_a+\epsi_{i-b+1}}$ for $a<i-b+1$, $(\Uu^i)_{a,i-a+1}=0$, and let $\Uu_s^i$ be the matrix obtained from $\Uu^i$ by deleting the $(i-s+1)$-th row and column. Then we can set
\begin{equation}
\begin{split}
a_s&=\sum\limits_{j=s+1}^ne_{\epsi_s-\epsi_j}\label{formula:Delta_i_D_n_net}e_{\epsi_s+\epsi_j},\\
\xi_{\epsi_{i-1}-\epsi_i}&=\sum\limits_{s=1}^ia_s\det\Uu_s^i.\\
\end{split}
\end{equation}
Finally, assume $m=n$ is even. In this case $\xi_{\epsi_{n-1}-\epsi_n}$ can be defined by
\begin{equation*}
\xi_{\epsi_{n-1}-\epsi_n}^2=\pm\begin{vmatrix}
e_{\epsi_1-\epsi_n}&e_{\epsi_1+\epsi_{n-1}}&\ldots&e_{\epsi_1+\epsi_2}&0\\
e_{\epsi_2-\epsi_n}&e_{\epsi_2+\epsi_{n-1}}&\ldots&0&-e_{\epsi_1+\epsi_2}\\
\vdots&\vdots&\iddots&\vdots&\vdots\\
e_{\epsi_{n-1}-\epsi_n}&0&\ldots&-e_{\epsi_2+\epsi_{n-1}}&-e_{\epsi_1+\epsi_{n-1}}\\
0&-e_{\epsi_{n-1}-\epsi_n}&\ldots&-e_{\epsi_2-\epsi_n}&-e_{\epsi_1-\epsi_n}\\
\end{vmatrix}
\end{equation*}
(our normalization is such that the term $e_{\epsi_{n-1}-\epsi_n}\xi_{\epsi_{n-3}+\epsi_{n-2}}$ enters $\xi_{\epsi_{n-1}+\epsi_n}$ with coefficient 1).

iv) $\Phi=B_n$. If $i\leq n-1$ is odd, then $\xi_{\epsi_i+\epsi_{i+1}}$ can be defined via formula~(\ref{formula:Delta_i_D_n_pf}). Next, for even $i$, we can define $\xi_{\epsi_i-\epsi_{i+1}}$ via formula (\ref{formula:Delta_i_D_n_net}) with $$b_s=\sum\nolimits_{j=s+1}^ne_{\epsi_s-\epsi_j}e_{\epsi_s+\epsi_j}+e_{\epsi_s}^2/4$$ instead of $a_s$. Finally, assume $n=m$ is odd. Then $\xi_{\epsi_n}$ can be defined by
\begin{equation*}
\xi_{\epsi_n}^2=\pm\begin{vmatrix}
e_{\epsi_1}&e_{\epsi_1+\epsi_n}&\ldots&e_{\epsi_1+\epsi_2}&0\\
e_{\epsi_2}&e_{\epsi_2+\epsi_n}&\ldots&0&-e_{\epsi_1+\epsi_2}\\
\vdots&\vdots&\iddots&\vdots&\vdots\\
e_{\epsi_n}&0&\ldots&-e_{\epsi_2+\epsi_n}&-e_{\epsi_1+\epsi_n}\\
0&-e_{\epsi_n}&\ldots&-e_{\epsi_2}&-e_{\epsi_1}\\
\end{vmatrix}
\end{equation*}
(our normalization is such that the term $e_{\epsi_n}\xi_{\epsi_{n-2}+\epsi_{n-1}}$ enters $\xi_{\epsi_n}$ with coefficient 1).

For $A_{n-1}$ and $C_n$, we denote $\Delta_{\beta}=\sigma(\xi_{\beta})$ for $\beta\in\Bu$. Since all $e_{\alpha}$ involved in $\xi_{\beta}$ (i.e., $e_{\alpha}$ which appear in a term of $\xi_{\beta}$) commute, we conclude that  $\Delta_{\beta}$, $\beta\in\Bu$, is defined as an element of  $U(\nt)$ again by formulas (\ref{formula:Delta_i_A_n}), (\ref{formula:Delta_i_C_n}) for $A_{n-1}$, $C_n$ respectively. For $B_n$ and $D_n$, we denote $P_{\beta}=\sigma(\xi_\beta)$ for $\beta\in\Bu$. If $i$ is odd, all $e_{\alpha}$ involved in $\xi_{\epsi_i+\epsi_{i+1}}$ commute, so the polynomial $P_{\epsi_i+\epsi_{i+1}}$ is defined as an element of $U(\nt)$ again by formula (\ref{formula:Delta_i_D_n_pf}).

\defi{For $A_{n-1}$ and $C_n$ (respectively, for $B_n$ and $D_n$), we call $\Delta_{\beta}$ (respectively, $P_{\beta}$), $\beta\in\Bu$, the \emph{canonical generators} of~$Z(\nt)$.}

\sst{Centrally generated ideals of $U(\nt)$} Throughout\label{sst:ideals_fd} this subsection $\gt$ and $\nt$ are as in Sub\-sec\-tion~\ref{sst:center_fd}. By definition, an ideal $J\subseteq U(\nt)$ is \emph{primitive} if $J$ is the annihilator of a simple $\nt$-module. Here we describe all primitive \emph{centrally generated} ideals of $U(\nt)$, i.e., all primitive ideals $J$ generated (as ideals) by their intersections $J\cap Z(\nt)$ with the center $Z(\nt)$ of $U(\nt)$.

In the 1960s A.~Kirillov, B.~Kostant and J.-M. Souriau discovered that the orbits of the coadjoint action play a crucial role in the representation theory of $B$ and $N$ (see, e.g., \cite{Kirillov1}, \cite{Kirillov2}). Works of J. Dixmier, M. Duflo, M. Vergne, O. Mathieu, N. Conze and R. Rentschler led to the result that the orbit method provides a nice description of primitive ideals of the universal enveloping algebra of a nilpotent Lie algebra (in particular, of $\nt$). Let us describe this in detail.

Let $\nt^*$ be the dual space of $\nt$. To any linear form $\lambda\in\nt^*$ one can assign a bilinear form $\beta_{\lambda}$ on $\nt$ by putting $\beta_{\lambda}(x,y)=\lambda([x,y])$. A subalgebra $\pt\subseteq\nt$ is a~\emph{polarization of $\nt$ at} $\lambda$ if it is a maximal $\beta_{\lambda}$-isotropic subspace. By \cite{Vergne}, such a subalgebra always exists. Let $\pt$ be a polarization of $\nt$ at $\lambda$, and $W$ be the one-dimensional representation of $\pt$ defined by $x\mapsto\lambda(x)$. Then the annihilator $J(\lambda)=\Ann_{U(\nt)}{V}$ of the induced representation $V = U(\nt) \otimes_{U(\pt)} W$ is a primitive two-sided ideal of $U(\nt)$. It turns out that $J(\lambda)$ depends only on $\lambda$ and not on the choice of polarization. Further, $J(\lambda)=J(\mu)$ if and only if the coadjoint $N$-orbits of $\lambda$ and $\mu$ coincide. Finally, the \emph{Dixmier map} $$\Du\colon\nt^*\to\Prim U(\nt),~\lambda\mapsto J(\lambda),$$ induces a homeomorphism between $\nt^*/N$ and $\Prim U(\nt)$, where the latter set is endowed with the Jacobson topology. (See \cite{Dixmier2}, \cite{Dixmier3}, \cite{BorhoGabrielRentschler} for the details.)

In addition, it is well known that the following conditions on an ideal $J\subset U(\nt)$ are equivalent\break \cite[Proposition 4.7.4, Theorem 4.7.9]{Dixmier3}:
\begin{equation}
\begin{split}
&\text{i) $J$ is \label{formula:conditions_Prim_fd}primitive;}\\
&\text{ii) $J$ is maximal;}\\
&\text{iii) the center of $U(\nt)/J$ is trivial;}\\
&\text{iv) $U(\nt)/J$ is isomorphic to a Weyl algebra of finitely many variables.}
\end{split}
\end{equation}
Recall that the Weyl algebra $\Au_s$ of $2s$ variables is the unital associative algebra with generators $p_i$,~$q_i$ for $1\leq i\leq s$, and relations $[p_i,q_i]=1$, $[p_i,q_j]=0$ for $i\neq j$, $[p_i,p_j]=[q_i,q_j]=0$ for all $i,~j$. Furthermore, in conditions (\ref{formula:conditions_Prim_fd}) we have $U(\nt)/J\cong\Au_s$ where $s$ equals one half of the dimension of the coadjoint $N$-orbit of $\lambda$, given that $J=J(\lambda)$.

Before formulating our first main result, we will present the classification of centrally generated primitive ideals of $U(\nt)$ for $A_{n-1}$ and $C_n$ from \cite{IgnatyevPenkov1}. Recall the definition of the Kostant cascade $\Bu$ from Subsection~\ref{sst:center_fd}. Let $\Delta\subset\Phi^+$ be the set of simple roots, i.e.,
\begin{equation*}
\Delta=\begin{cases}\{\epsi_1-\epsi_2,~\ldots,~\epsi_{n-1}-\epsi_n\}&\text{for }\Phi=A_{n-1},\\
\{\epsi_1-\epsi_2,~\ldots,~\epsi_{n-1}-\epsi_n,~\epsi_n\}&\text{for }\Phi=B_n,\\
\{\epsi_1-\epsi_2,~\ldots,~\epsi_{n-1}-\epsi_n,~2\epsi_n\}&\text{for }\Phi=C_n,\\
\{\epsi_1-\epsi_2,~\ldots,~\epsi_{n-1}-\epsi_n,~\epsi_{n-1}+\epsi_n\}&\text{for }\Phi=C_n.
\end{cases}
\end{equation*}
Note also that
\begin{equation*}
\Bu\setminus\Delta=\begin{cases}
\bigcup_{1\leq i\leq(n-1)/2}\{\epsi_{i}-\epsi_{n-i+1}\}&\text{for }\Phi=A_{n-1},\\
\bigcup_{1\leq i\leq n/2}\{\epsi_{2i-1}+\epsi_{2i+1}\}&\text{for }\Phi=B_n,\\
\bigcup_{1\leq i<n/2}\{\epsi_{2i-1}+\epsi_{2i+1}\}&\text{for }\Phi=D_n,\\
\Bu\setminus\{2\epsi_n\}&\text{for }\Phi=C_n.
\end{cases}
\end{equation*}.

Let $\{e_{\alpha}^*,~\alpha\in\Phi^+\}$ be the basis of $\nt^*$ dual to the basis $\{e_{\alpha},~\alpha\in\Phi^+\}$ of $\nt$.

\defi{To a map $\xi\colon\Bu\to\Cp$ we assign the linear form $f_{\xi}=\sum_{\beta\in\Bu}\xi(\beta)e_{\beta}^*\in\nt^*$. We call a form $f_{\xi}$ a \emph{Kostant form} if $\xi(\beta)\neq0$ for any $\beta\in\Bu\setminus\Delta$.}

Let $V$ be a simple $\nt$-module and $J=\Ann_{U(\nt)}{V}$ be the corresponding primitive ideal of $U(\nt)$. By a version of Schur's Lemma \cite{Dixmier4}, each central element of $U(\nt)$ acts on $V$ as a scalar operator. For $A_{n-1}$ and $C_n$ (respectively, for $B_n$ and $D_n$), let $c_{\beta}$ be the scalar corresponding to~$\Delta_{\beta}$ (respectively, to~$P_{\beta}$) for $\beta\in\Bu$, and $J_c$ be the ideal of $U(\nt)$ generated by all $\Delta_{\beta}-c_\beta$ (respectively, by all $P_{\beta}-c_{\beta}$), $\beta\in\Bu$. Clearly, $J_c\subseteq J$. Further, since $Z(\nt)$ is a polynomial ring and the center of $U(\nt)/J$ is trivial, $J$~is centrally generated if and only if $J=J_c$.

The following result was proved in \cite[Theorem 3.1]{IgnatyevPenkov1}.\newpage

\mtheo{Suppose $\Phi$ is of type $A_{n-1}$ or $C_n$. The following conditions on a primitive ideal $J\subset U(\nt)$ are equivalent\textup{:}
\begin{equation*}
\begin{split}
&\text{\textup{i)} $J$ is centrally generated \textup{(}or\textup{,} equivalently\textup{,} $J=J_c$\textup{)}};\\
&\text{\textup{ii)} the scalars $c_{\beta}$\textup, $\beta\in\Bu\setminus\Delta$\textup, are nonzero};\\
&\text{\textup{iii)} $J=J(f_{\xi})$ for a Kostant form $f_{\xi}\in\nt^*$}.\\
\end{split}
\end{equation*}
If these\label{theo:ideal_fd_An_Cn} conditions are satisfied\textup{,} then the map $\xi$ is reconstructed by $J$\textup{:}
\begin{equation}
\xi(\beta_k)=(-1)^{k+1}\label{formula_xi_via_c_fd_An_Cn}c_kc_{k-1}^{-1},
\end{equation}
where $c_0=1$ and $c_k=c_{\beta_k}$ for $k\geq 1$. Here $\beta_k=\epsi_k-\epsi_{n-k+1}$ for $\Phi=A_{n-1}$\textup{,} and  $\beta_k=2\epsi_k$ for $\Phi=C_n$.}


One of the key ingredients in the proof of this result is to check that if condition (ii) is satisfied then $J_c$ is primitive. To do this, in \cite{IgnatyevPenkov1} an explicit set of generators of the quotient algebra $U(\nt)/J_c$ was constructed. It turns out that these generators satisfy (up to scalars) the defining relations of the Weyl algebra $\Au_s$ for $s=|\Phi^+\setminus\Bu|/2$. Since $\Au_s$ is simple and, as one can easily check, $J_c\neq U(\nt)$, we conclude that $U(\nt)/J_c\cong\Au_s$, and, consequently, $J_c$ is primitive. Below we briefly recall the explicit formulas for these generators from \cite[Subsection 3.1]{IgnatyevPenkov1}.

We define the maps $\row\colon\Phi^+\to\Zp$ and $\col\colon\Phi^+\to\Zp$ by putting $\row(\epsi_i-\epsi_j)=\row(\epsi_i+\epsi_j)=\row(2\epsi_i)=\row(\epsi_i)=i$, $\col(\epsi_i+\epsi_j)=\col(2\epsi_j)=-j$. Let $\Ro_i=\{\alpha\in\Phi^+\mid\row(\alpha)=i\}$. For $\alpha\in\Phi^+$, where $\Phi=A_{n-1}$ or $C_n$, set
\begin{equation*}
\begin{split}
A(\alpha)&=\begin{cases}\bigcup\nolimits_{j+1\leq k\leq n-i+1}\{\epsi_j-\epsi_k\},&\text{if $\Phi=A_{n-1}$, $\alpha=\epsi_i-\epsi_j$, $j<n-i+1$},\\
\bigcup\nolimits_{n-j+1\leq k\leq i-1}\{\epsi_k-\epsi_i\},&\text{if $\Phi=A_{n-1}$, $\alpha=\epsi_i-\epsi_j$, $j>n-i+1$},\\
\bigcup\nolimits_{i\leq k\leq j-1}\{\epsi_k+\epsi_j\}\cup\Ro_j,&\text{if $\Phi=C_n$, $\alpha=\epsi_i-\epsi_j$},\\
\bigcup\nolimits_{i\leq k\leq j-1}\{\epsi_k-\epsi_j\},&\text{if $\Phi=C_n$, $\alpha=\epsi_i+\epsi_j$},\\
\end{cases}\\
\Bu(\alpha)&=\{\alpha\}\cup\{\beta\in\Bu\mid\row(\beta)<\row(\alpha)\},\\
R(\alpha)&=\{\row(\gamma),~\gamma\in\Bu(\alpha)\},~C(\alpha)=\{\col(\gamma),~\gamma\in\Bu(\alpha)\}.
\end{split}
\end{equation*}
Define a matrix $\Uu$ with entries from $U(\nt)$ by the following rule.

\begin{center}
\begin{tabular}{|l|l|l|}
\hline
$\Phi$&Size of $\Uu$&$\Uu$\\
\hline\hline
$A_{n-1}$&$n\times n$&$\Uu_{i,j}=e_{\epsi_i-\epsi_j}$ for $1\leq i<j\leq n$,\\
&&$\Uu_{i,j}=0$ otherwise\\
\hline
$C_n$&$2n\times2n$&$\Uu_{i,j}=-\Uu_{-j,-i}=e_{\epsi_i-\epsi_j}$, $\Uu_{i,-j}=\Uu_{j, -i}=e_{\epsi_i+\epsi_j}$ for $1\leq i<j\leq n$,\\
&&$\Uu_{i,-i}=2e_{\epsi_{2i}}$, $1\leq i\leq n$, $\Uu_{i,j}=0$ otherwise\\
\hline
\end{tabular}
\end{center}
As before, for $C_n$, we index the rows and the columns by the numbers $1,~\ldots,~n,~-n,~\ldots,~-1$.

Denote by $\Delta_{\alpha}$ the element of $U(\nt)$, which equals the minor of $\Uu$ with rows $R(\alpha)$ and columns $C(\alpha)$. Note that the variables involved in each $\Delta_{\alpha}$ commute. For example, let $\Phi=A_{n-1}$, $n=8$, $\alpha=\epsi_3-\epsi_4$. On the picture below $\alpha$ is marked by $\bullet$, the roots from $\Bu$ are marked by $\otimes$'s, and the roots $\gamma$ such that $e_{\gamma}$ is involved in $\Delta_{\alpha}$ are grey:
\begin{equation*}\predisplaypenalty=0
\mymatrix{
\pho& \Lft{2pt}\Bot{2pt}\pho& \pho& \gray\pho& \pho& \pho& \gray\pho& \gray\otimes\\
\pho& \pho& \Lft{2pt}\Bot{2pt}\pho& \gray\pho& \pho& \pho& \gray\otimes& \gray\pho\\
\pho& \pho& \pho& \Lft{2pt}\Bot{2pt}\gray\bullet& \pho& \otimes& \gray\pho& \gray\pho\\
\pho& \pho& \pho& \pho& \Lft{2pt}\Bot{2pt}\otimes& \pho& \pho& \pho\\
\pho& \pho& \pho& \pho& \pho& \Lft{2pt}\Bot{2pt}\pho& \pho& \pho\\
\pho& \pho& \pho& \pho& \pho& \pho& \Lft{2pt}\Bot{2pt}\pho& \pho\\
\pho& \pho& \pho& \pho& \pho& \pho& \pho& \Lft{2pt}\Bot{2pt}\pho\\
\pho& \pho& \pho& \pho& \pho& \pho& \pho& \pho\\
}\qquad
\Delta_{\epsi_3-\epsi_4}=\begin{vmatrix}
e_{\epsi_1-\epsi_4}&e_{\epsi_1-\epsi_7}&e_{\epsi_1-\epsi_8}\\
e_{\epsi_2-\epsi_4}&e_{\epsi_2-\epsi_7}&e_{\epsi_2-\epsi_8}\\
e_{\epsi_3-\epsi_4}&e_{\epsi_3-\epsi_7}&e_{\epsi_3-\epsi_8}\\
\end{vmatrix}.\postdisplaypenalty=10000\end{equation*}
Note also that for $\alpha\in\Bu$ this definition agrees with the definition given in Subsection~\ref{sst:center_fd}.

Next, let $\wh a$ be the image of an element $a\in U(\nt)$ under the canonical projection $U(\nt)\sur U(\nt)/J_c$. It turns out that the elements $\wh\Delta_{\alpha}$, $\alpha\in\Phi^+\setminus\Bu$, generate the quotient algebra $U(\nt)/J_c$. Furthermore, given $\alpha\in\Phi^+$ with $\row(\alpha)=i$, set
$p_{\alpha}=\wh\Delta_{\alpha},~q_{\alpha}=(-1)^{i+1}c_i^{-1}c_{i-1}^{-1}\wh\Delta_{\beta_i-\alpha}.$  Assume that $\Phi=A_{n-1}$ and $\col(\alpha)<n-i+1$, $\col(\gamma)<n-i+1$, or that $\Phi=C_n$ and $\col(\alpha)>0$, $\col(\gamma)>0$. Then, as it is shown in the proof of \cite[Lemma 1.5]{IgnatyevPenkov1},
\begin{equation*}
[p_{\alpha},q_{\gamma}]=\begin{cases}
1,&\text{if }\alpha=\gamma,\\
0&\text{otherwise}.\\
\end{cases}
\end{equation*}
Thus, $U(\nt)/J_c\cong\Au_s$ for $s=|\Phi^+\setminus\Bu|/2$, as required.

Now we will formulate our first main result (cf. Theorem~\ref{theo:ideal_fd_An_Cn}).
\mtheo{Suppose $\Phi$ is of type $B_n$ or $D_n$. The following conditions on a primitive ideal $J\subset U(\nt)$ are equivalent\textup{:}
\begin{equation*}
\begin{split}
&\text{\textup{i)} $J$ is centrally generated \textup{(}or\textup{,} equivalently\textup{,} $J=J_c$\textup{)}};\\
&\text{\textup{ii)} the scalars $c_{\beta}$\textup, $\beta\in\Bu\setminus\Delta$\textup, are nonzero};\\
&\text{\textup{iii)} $J=J(f_{\xi})$ for a Kostant form $f_{\xi}\in\nt^*$}.\\
\end{split}
\end{equation*}
If these\label{theo:ideal_fd_Bn_Dn} conditions are satisfied\textup{,} then the map $\xi$ is reconstructed by $J$\textup{:}
\begin{equation}
\xi(\beta)\label{formula_xi_via_c_fd_Bn_Dn}=
\begin{cases}
c_{\epsi_1+\epsi_2}&\text{for }\beta=\epsi_1+\epsi_2,\\
c_{\epsi_{2k-1}+\epsi_{2k}}c_{\epsi_{2k-3}+\epsi_{2k-2}}^{-1}&\text{for }\beta=\epsi_{2k-1}+\epsi_{2k},~k\geq2,\\
c_{\epsi_1-\epsi_2}c_{\epsi_1+\epsi_2}^{-1}&\text{for }\beta=\epsi_1-\epsi_2,\\
c_{\epsi_{2k-1}-\epsi_{2k}}c_{\epsi_{2k-1}+\epsi_{2k}}^{-1}c_{\epsi_{2k-3}+\epsi_{2k-2}}^{-1}&\text{for }\beta=\epsi_{2k-1}-\epsi_{2k},~k\geq2,~2k<n\text{ if }\Phi=D_n,\\
c_{\epsi_{n-1}-\epsi_n}c_{\epsi_{n-1}+\epsi_n}^{-1}&\text{for }\beta=\epsi_1-\epsi_n,~\Phi=D_n,~n\text{ even},\\
c_{\epsi_n}c_{\epsi_{n-2}+\epsi_{n-1}}^{-1}&\text{for }\beta=\epsi_n,~\Phi=B_n,~n\text{ odd}.\\
\end{cases}
\end{equation}}

Again, as for $A_{n-1}$ and $C_n$, one of the most important steps is to prove that if condition (ii) is satisfied then the centrally generated ideal $J_c$ is primitive. To do this in the orthogonal case, we need some additional technical (but important) facts. Put $\beta_1=\epsi_1+\epsi_2$, $\beta_2=\epsi_1-\epsi_2$, $e=e_{\beta_2}$, and denote $\wt\nt=\langle e_{\alpha},~\alpha\in\Phi^+$, $\row(\alpha)>2\rangle_{\Cp},~\kt=\langle e_{\alpha},~\row(\alpha)\leq2,~\alpha\neq\beta_2\rangle_{\Cp}$. For $\Phi=B_n$ (res\-pec\-ti\-vely, for $\Phi=D_n$), $\wt\nt$ is isomorphic to the nilradical of a Borel subalgebra of the simple Lie algebra with the root system $\wt\Phi=B_{n-2}$ (res\-pec\-tively, $\wt\Phi=D_{n-2}$), and $\kt$ is isomorphic to the $(4n-5)$-dimen\-sional (respectively, $(4n-7)$-dimensional) Heisenberg algebra. Indeed, $[e_{\epsi_1-\epsi_i},e_{\epsi_2+\epsi_i}]=e_{\beta_1}$ for $3\leq i\leq n$, $[e_{\epsi_1},e_{\epsi_2}]=2e_{\beta_1}$, and $[e_{\alpha},e_{\gamma}]=0$ for all other $\alpha,~\gamma\in(\Ro_1\cup\Ro_2)\setminus\{\beta_2\}$.

Given $c_1\in\Cp^{\times}$, denote by $J_1$ the ideal of $U(\kt)$ generated by $e_{\beta_1}-c_1$, then, clearly, $U(\kt)/J_1\cong\Au_s$, where $s=2n-3$ for $B_n$ and $s=2n-4$ for $D_n$. Denote $\wt\nt_e=\wt\nt\oplus\Cp e$. Note that $\kt$ is an ideal of the Lie algebra $\nt$, so, given $x\in\wt\nt_e$, one can consider $\ad{x}$ as a derivation of $\kt$. Since $e_{\beta_1}-c_1$ is a central element of $U(\nt)$, one has $\ad{x}(J_1)\subseteq J_1$, so $\ad{x}$ can be considered as a derivation of $\Au_s$. It is well known (see, e.g., \cite[Lemma~10.1.2, Lemma 10.1.3]{Dixmier3}) that there exist the unique element $\theta(x)\in\Au_s$ such that $\ad{x}(y)=[\theta(x),y]$ for all $y\in\Au_s$, and $\theta\colon\wt\nt_e\to\Au_s$ is a morphism of Lie algebras. Furthermore, there exists the unique epimorphism of associative algebras $$r\colon U(\nt)\sur U(\wt\nt_e)\otimes\Au_s$$ such that $r(y)=1\otimes \overline y$ for $y\in\kt$ and $r(x)=x\otimes1+1\otimes\theta(x)$ for $x\in\wt\nt_e$. (Here $\overline a$ is the image of an element $a\in U(\kt)$ under the canonical projection $U(\kt)\sur U(\kt)/J_1\cong\Au_s$.) It turns out that the kernel of $r$ coincides with the ideal $J_0$ of $U(\nt)$ generated by $e_{\beta_1}-c_1$ \cite[Lemma 10.1.5]{Dixmier3}. One can easily check by direct computation that $r(P_{\beta_2})=c_1e\otimes1$. We are now ready to formulate and prove the main technical proposition needed for the proof of Theorem~\ref{theo:ideal_fd_Bn_Dn}.\newpage

\propp{Let\textup, as\label{prop:J_c_primitive_Bn_Dn} above\textup, $J_c$ be the ideal of $U(\nt)$ generated by ${P}_{\beta}-c_{\beta}$\textup, $\beta\in\Bu$\textup, with $c_{\beta}\neq0$ for $\beta\in\Bu\setminus\Delta$. Then $J_c$ is primitive.}{Put $c_1=c_{\beta_1}$. Since $r$ is surjective, $r(J_c)$ is an ideal of $U(\wt\nt_e)\otimes\Au_s$ generated by $r({P}_{\beta})-c_{\beta}$, $\beta\in\Bu\setminus\{\beta_1\}$, or, equivalently, by $e\otimes1-c_{\beta_1}^{-1}c_{\beta_2}$ and by $r({P}_{\beta})-c_{\beta}$, $\beta\in\wt\Bu=\Bu\setminus\{\beta_1,~\beta_2\}$. Note that $\wt\Bu$ is the Kostant cascade of $\wt\Phi$ (under the natural identification of $\wt\Phi$ with the root subsystem $\pm\{\alpha\in\Phi^+\mid\row(\alpha)>2\}$ of $\Phi$). Denote by $\wt{P}_{\beta}$, $\beta\in\wt\Phi$, the set of canonical generators of $Z(\wt\nt)$. Our first goal is to check that, up to nonzero scalar, $r({P}_{\beta})$ coincides with $\wt{P}_{\beta}\otimes1$ for all $\beta\in\wt\Bu$.

To check this fact, we will use (\ref{formula:weights_of_xi_i}). Pick a root $\beta\in\wt\Bu$. Since $r$ is surjective, $r({P}_{\beta})$ is a central element of $U(\wt\nt_e)\otimes\Au_s$. The center of this algebra has the form $Z(\wt\nt_e)\otimes\Cp$, so in fact $r({P}_{\beta})\in Z(\wt\nt_e)\otimes\Cp$. The Kostant's description of the center of the enveloping algebra is valid for the nilradical of a Borel subalgebra of any semisimple Lie algebra with the Kostant cascade being the union of the Kostant cascades of the simple summands \cite{Kostant1}, \cite{Kostant2}. Note that $\wt\nt_e$ is isomorphic to the nilradical of a Borel subalgebra of the semisimple Lie algebra $\wt\gt_e$ with the root system $\wt\Phi\times A_1$, and the Kostant cascade of $\wt\nt_e$ is $\wt\Bu\cup\{\beta_2\}$. Denote by $\wt\htt_e$ the corresponding Cartan subalgebra of $\wt\gt_e$. By \cite[Theorem 6]{Kostant2} (see also \cite[Lemma 4.4]{Joseph1}), $Z(\nt)$ (respectively, $Z(\wt\nt_e)$) is a direct sum of 1-dimensional weight spaces of~$\htt$ (respectively, of $\wt\htt_e$) with respect to the adjoint action of the corresponding Cartan subalgebras. Since $[h,e_{\beta_1}]=0$ for all $h\in\wt\htt_e$, the algebra $\wt\htt_e$ naturally acts on $\Au_s$, and so on $U(\wt\nt_e)\otimes\Au_s$. We define the result of this action by $h.x$, $h\in\wt\htt_e$, $x\in U(\wt\nt_e)\otimes\Au_s$. Hence, according to (\ref{formula:weights_of_xi_i}), it is enough to check that $r({P}_{\beta})$ is a nonzero $\wt\htt_e$-weight element of weight $\wt\vpi_{\beta}$, where
\begin{equation*}
\wt\vpi_{\beta}=\begin{cases}
\epsi_3+\ldots+\epsi_{2k}&\text{for }\beta=\epsi_{2k-1}+\epsi_{2k},~k\geq2,\\
\epsi_3+\ldots+\epsi_n&\text{for }\beta=\epsi_n,~\Phi=B_n,~n\text{ odd},\\
\epsi_3+\ldots+\epsi_{n-1}-\epsi_n&\text{for }\beta=\epsi_{n-1}-\epsi_n,~\Phi=D_n,~n\text{ even},\\
2\epsi_3+\ldots+2\epsi_{2k-1}&\text{for }\beta=\epsi_{2k-1}-\epsi_{2k},~k\geq2,~2k<n\text{ if }\Phi=D_n.\\
\end{cases}
\end{equation*}

To prove that $r({P}_{\beta})$ is an $\wt\htt_e$-weight element of weight $\wt\vpi_{\beta}$, denote the result of the natural (adjoint) action of $\wt\htt_e$ on $U(\nt)$ by $h\cdot x$, $h\in\wt\htt_e$, $x\in U(\nt)$. As above, since $\wt\htt_e\cdot e_{\beta_1}=0$, the algebra $\wt\htt_e$ naturally acts on $U(\nt)/J_0$ by the formula $h\cdot r(x)=r(h\cdot x)$. We claim that this action coincides with the action of $\wt\htt_e$ on $U(\wt\nt_e)\otimes\Au_s$ defined above, i.e., that $h.x=h\cdot x$ for all $h\in\wt\htt_e$, $x\in U(\wt\nt_e)\otimes\Au_s$. Indeed, if $y\in\kt$, then $h\cdot r(y)=r([h,y])=1\otimes\overline{[h,y]}=h.(1\otimes y)=h.r(y)$, as required. On the other hand, if $e_{\alpha}\in\wt\nt_e$ for some $\alpha\in\Phi^+$, then, by \cite[Subsection 4.8]{Joseph1}, $\theta(e_{\alpha})$ is a linear combination of elements of the form $\overline{e}_{\alpha+\gamma}\overline{e}_{\beta_1-\gamma}$, $\gamma\in(\Ro_1\cup\Ro_2)\setminus\{\beta_1,\beta_2\}$. We conclude that $h(1\otimes\theta(e_{\alpha}))=(\beta_1+\alpha)(h)1\otimes\theta(\alpha)=\alpha(h)1\otimes\theta(e_{\alpha})$, because $\beta_1(\wt\htt_e)=0$. Thus,
$$h.r(e_{\alpha})=[h,e_{\alpha}]\otimes1+1\otimes\alpha(h)\theta(e_{\alpha})=\alpha(h)1\otimes r(e_{\alpha})=h\cdot r(e_{\alpha}).$$ It remains to note that ${P}_{\beta}$ is an $\htt$-weight element of $U(\nt)$ of weight $\beta_1+\wt\vpi_{\beta}$, but $\beta_1(\wt\htt_e)=0$. To show that $r({P}_{\beta})\neq0$, recall that the kernel of $r$ is $J_0$. If ${P}_{\beta}\in J_0$ (i.e., if ${P}_{\beta}=(e_{\beta_1}-c_{\beta_1})a$ for some $a\in U(\nt)$), then, clearly, $a\in Z(\nt)$. But this contradicts the fact that ${P}_{\beta}$ and $e_{\beta_1}$ are algebraically independent.

So, given $\beta\in\wt\Bu$, there exists the unique $a_{\beta}\in\Cp^{\times}$ such that $r({P}_{\beta})=a_{\beta}\wt{P}_{\beta}\otimes1$. Consequently, $r(J_c)$ is generated by $e\otimes1-c_{\beta_1}^{-1}c_{\beta_2}$ and by $\wt{P}_{\beta}\otimes1-\wt c_{\beta}$, $\beta\in\wt\Bu$, where $\wt c_{\beta}=a_{\beta}^{-1}c_{\beta}$. In particular, $\wt c_{\beta}\neq0$ if $\beta\in\wt\Bu$ is not a simple root of $\wt\Phi^+$. Now we will use the induction on $\rk\Phi$ to prove that $J_c$ is primitive. One can check the base $\rk\Phi\leq3$ directly using explicit formulas from Subsection~\ref{sst:center_fd}. Denote by $\wt J_c$ the ideal of $U(\wt\nt)$ generated by $\wt{P}_{\beta}-\wt c_{\beta}$, $\beta\in\wt\Bu$. By the inductive assumption, $\wt J_c$ is a primitive ideal of $U(\wt\nt)$, so $U(\wt\nt)/\wt J_c\cong\Au_t$ for certain $t$. Taking in account that $U(\wt\nt_e)=U(\wt\nt)\otimes\Cp[e]$, and that the quotient algebra of $\Cp[e]$ modulo the ideal generated by $e-c_{\beta_1}^{-1}c_{\beta_2}$ is isomorphic to $\Cp$, we conclude that
\begin{equation*}
U(\nt)/J_c\cong(U(\nt)/J_0)/r(J_c)\cong(U(\wt\nt_e)\otimes\Au_s)/r(J_c)=(U(\wt\nt)/\wt J_c)\otimes\Cp\otimes\Au_s\cong\Au_t\otimes\Au_s\cong\Au_{t+s}.
\end{equation*}
Thus, $J_c$ is primitive. The proof is complete.}\newpage

Now, using this proposition we will prove our first main result, Theorem~\ref{theo:ideal_fd_Bn_Dn}. Note that each element of $S(\nt)$ can be considered as a polynomial function on $\nt^*$ via the natural isomorphism $(\nt^*)^*\cong\nt$.

\medskip\textsc{Proof of Theorem~\ref{theo:ideal_fd_Bn_Dn}}. $\mathrm{(ii)}\Longrightarrow\mathrm{(iii)}$. Define $\xi$ be formula (\ref{formula_xi_via_c_fd_Bn_Dn}), and set $f_{\xi}$ to be the corresponding Kostant form. By \cite[6.6.9 (c)]{Dixmier3}, $J(f_{\xi})$ contains ${P}_{\beta}-\xi_{\beta}(f_{\xi})$ for all $\beta\in\Bu$. Using explicit formulas for $\xi_{\beta}$ from Subsection~\ref{sst:center_fd}, one can immediately check that $\xi_{\beta}(f_{\xi})=c_{\beta}$ for all $\beta\in\Bu$. Thus, both $J$ and $J(f_{\xi})$ contain the centrally generated ideal $J_c$. But, thanks to Proposition~\ref{prop:J_c_primitive_Bn_Dn}, $J_c$ is primitive (and, consequently, maximal).We conclude that $J=J_c=J(f_{\xi})$, as required.

$\mathrm{(iii)}\Longrightarrow\mathrm{(i)}$. Again by \cite[6.6.9 (c)]{Dixmier3}, $c_{\beta}=\xi_{\beta}(f_{\xi})$ for all $\beta\in\Bu$, hence $c_{\beta}$'s satisfy (\ref{formula_xi_via_c_fd_Bn_Dn}). Both $J$~and~$J(f_{\xi})$ contain the centrally generated ideal $J_c$, and condition (ii) is satisfied, so $J_c$ is primitive and $J=J_c=J(f_{\xi})$, as required.

$\mathrm{(i)}\Longrightarrow\mathrm{(ii)}$. Assume, to the contrary, that some scalars $c_{\beta}$, $\beta\in\Bu\setminus\Delta$, equal zero. Suppose that $i_1$ is the minimal number such that $c_{\beta}=0$ for $\beta\in(\Bu\setminus\Delta)\cap\Ro_{i_1}$. Now, as in the proof of \cite[Theorem~3.1]{IgnatyevPenkov1}, define inductively two (finite) sequences $\{i_j\}$ and $\{k_j\}$ of positive integers by the following rule. If $i_j$ is already defined and there exists $k>i_j$ such that $c_{\beta}\neq0$ for $\beta\in\Bu\cap\Ro_k$ with $\col(\beta)\leq0$, then set $k_j$ to be the minimal among all such $k$. Similarly, if $k_j$ is already defined and there exists $i>k_j$ such that $c_{\beta}=0$ for $\beta\in(\Bu\setminus\Delta)\cap\Ro_i$, then set $i_{j+1}$ to be the minimal among all such $i$.

To each $j$ for which both $i_j$ and $k_j$ exist we assign the set of roots
\begin{equation*}
\begin{split}
\Gamma_j&=\{\epsi_{i_j}+\epsi_{i_j+2},~\epsi_{i_j}-\epsi_{i_j+2}\}\cup\bigcup\limits_{i_j<s<k_j-1,~s\text{ even}}\{\epsi_s+\epsi_{s+3},~\epsi_s-\epsi_{s+3}\}\cup\Gamma_j',\text{ where}\\
\Gamma_j'&=\begin{cases}
\{\epsi_{n-1}\},&\text{if $\Phi=B_n$ and $k_j=n$},\\
\{\epsi_{k_j-1}+\epsi_{k_j+1},~\epsi_{k_j-1}-\epsi_{k_j+1}\}&\text{otherwise}.
\end{cases}
\end{split}
\end{equation*}
Denote the length of the sequence $\{i_j\}$ by $l_I$, then the length $l_K$ of the sequence $\{k_J\}$ is either $l_I$ or $l_I-1$. If $i_j$ exists and $k_j$ does not exist (i.e., if $l_K=l_I-1$ and $j=l_I$), then put
\begin{equation*}
\begin{split}
\Gamma_j&=\{\epsi_{i_j}+\epsi_{i_j+2},~\epsi_{i_j}-\epsi_{i_j+2}\}\cup\bigcup\limits_{i_j<s<n-2,~s\text{ even}}\{\epsi_s+\epsi_{s+3},~\epsi_s-\epsi_{s+3}\}\cup\Gamma_j',\text{ where}\\
\Gamma_j'&=\begin{cases}
\{\epsi_{n-2}\},&\text{if $\Phi=B_n$, $n$ even},\\
\{\epsi_{n-2}-\epsi_n,~\epsi_{k_j-1}-\epsi_{k_j+1}\},&\text{if $\Phi=D_n$, $n$ even},\\
\varnothing&\text{otherwise}.
\end{cases}
\end{split}
\end{equation*}
Finally, put
\begin{equation*}
X=\left(\Bu\cup\bigcup_{j=1}^{l_I}\Gamma_j\right)\setminus\left(\bigcup\{\Ro_i\mid c_{\beta}=0\text{ for }\beta\in(\Bu\setminus\Delta)\cap\Ro_i\}\cup\bigcup_{j=1}^{l_K}\Ro_{k_j}\right).
\end{equation*}

Given a map $\vfi\colon X\to\Cp$, we denote $$\mu_{\vfi}=\sum_{\alpha\in X}\vfi(\alpha)e_{\alpha}^*.$$ By \cite[6.6.9 (c)]{Dixmier3}, ${P}_{\beta}-\xi_{\beta}(\mu_{\vfi})\in J(\mu_{\vfi})$, $\beta\in\Bu$. Using explicit formulas for ${P}_{\beta}$'s from Subsection~\ref{sst:center_fd}, one can easily construct a map $\vfi_1\colon X\to\Cp$ such that $\vfi_1(\beta)\neq0$ if $\col(\beta)\leq0$, and $\xi_{\beta}(\mu_{\vfi_1})=c_{\beta}$ for all $\beta\in\Bu$. Let $\vfi_2\colon X\to\Cp$ be the map for which $\vfi_2(\epsi_{i_1}+\epsi_{i_1+2})=-\vfi_1(\epsi_{i_1}+\epsi_{i_1+2})$ and $\vfi_2(\beta)=\vfi_1(\beta)$ for all other $\beta\in\Bu$. One can check that $\xi_{\beta}(\mu_{\vfi_1})=\xi_{\beta}(\mu_{\vfi_2})$ for all $\beta\in\Bu$. It follows from \cite[Theorem 1.4]{Panov2} that the coadjoint orbits of $\mu_{\vfi_1}$ and $\mu_{\vfi_1}$ are disjoint, so $J(\mu_{\vfi_1})\neq J(\mu_{\vfi_2})$. On the other hand, both $J(\mu_{\vfi_1})$ and $J(\mu_{\vfi_2})$ contain $J=J_c$, and this contradicts the maximality of~$J$. The equivalence of (i), (ii), (iii) is now proved. The fact that the map $\xi$ is reconstructed by $J$ via formula~(\ref{formula_xi_via_c_fd_An_Cn}) follows from the proof of the implication $\mathrm{(iii)}\Longrightarrow\mathrm{(i)}$.\hfill$\square$

\medskip In the rest of this subsection, given a primitive ideal $J$ of $U(\nt)$, we will construct an explicit set of generators of the quotient algebra $U(\nt)$ satisfying the defining relations of the Weyl algebra of corresponding rank. Pick $\alpha\in\Phi^+$ and put
\begin{equation*}
\begin{split}
A(\alpha)&=\begin{cases}
\Ro_j\cup\bigcup\nolimits_{i<k<j}\{\epsi_k+\epsi_j\}\cup\{\epsi_{i-1}+\epsi_j,~\epsi_{i-1}-\epsi_i\}&\text{for }\alpha=\epsi_i-\epsi_j,~i<j,~i\text{ even},\\
\Ro_j\cup\bigcup\nolimits_{i<k<n}\{\epsi_k+\epsi_j\}&\text{for }\alpha=\epsi_i-\epsi_j,~i<j,~i\text{ odd},\\
\bigcup\nolimits_{i<k<j}\{\epsi_k-\epsi_j\}\cup\{\epsi_{i-1}-\epsi_j,~\epsi_{i-1}-\epsi_i\}&\text{for }\alpha=\epsi_i+\epsi_j,~i<j,~i\text{ even},\\
\bigcup\nolimits_{i<k<j}\{\epsi_k-\epsi_j\}\cup\{\epsi_{i-1}-\epsi_j\}&\text{for }\alpha=\epsi_i+\epsi_j,~i<j,~i\text{ odd},\\
\bigcup\nolimits_{i<k\leq n}\{\epsi_k\}&\text{for }\alpha=\epsi_k,
\end{cases}\\
\Bu(\alpha)&=\{\gamma\in\Bu\setminus\Delta\mid\row(\gamma)\leq\row(\alpha)\},~R(\alpha)=\{\row(\gamma)\label{formula:A_B_R_C_Bn_Dn_fd},~-\col(\gamma),~\gamma\in\Bu(\alpha)\},\\
C(\alpha)&=\begin{cases}
-R(\alpha)\cup\{\col(\alpha)\}\setminus\{-\row(\alpha)+1\},&\text{if }\row(\alpha)\text{ is even},\\
-R(\alpha)\cup\{\col(\alpha)\}\setminus\{-\row(\alpha)-1\},&\text{if }\row(\alpha)\text{ is odd}.\\
\end{cases}
\end{split}
\end{equation*}
Define the matrix $\Uu$ by the following rule.
\begin{center}
\begin{tabular}{|l|l|l|}
\hline
$\Phi$&Size of $\Uu$&$\Uu$\\
\hline\hline
$B_n$&$(2n+1)\times(2n+1)$& $\Uu_{i,j}=-\Uu_{-j,-i}=e_{\epsi_i-\epsi_j}$, $\Uu_{i,-j}=-\Uu_{j, -i}=e_{\epsi_i+\epsi_j}$, $1\leq i<j\leq n$,\\
&&$\Uu_{i,0}=-\Uu_{0,i}=e_{\epsi_i}$, $1\leq i\leq n$, $\Uu_{i,j}=0$ otherwise\\
\hline
$D_n$&$2n\times2n$& $\Uu_{i,j}=-\Uu_{-j,-i}=e_{\epsi_i-\epsi_j}$, $\Uu_{i,-j}=-\Uu_{j, -i}=e_{\epsi_i+\epsi_j}$, $1\leq i<j\leq n$,\\
&&$\Uu_{i,j}=0$ otherwise\\
\hline
\end{tabular}
\end{center}
Denote by $\Delta_{\alpha}$ the element of $U(\nt)$, which equals the minor of the matrix $\Uu$ with the set of rows $R(\alpha)$ and the set of columns $C(\alpha)$. Note that variables in each summand of each $\Delta_{\alpha}$ commute. For instance, if $\Phi=D_n$, $n\geq6$, and $\alpha=\epsi_4-\epsi_6$, then
\begin{equation*}
\Delta_{\alpha}=\begin{vmatrix}
e_{\epsi_1-\epsi_6}&e_{\epsi_1+\epsi_4}&e_{\epsi_1+\epsi_2}&0\\
e_{\epsi_2-\epsi_6}&e_{\epsi_2+\epsi_4}&0&-e_{\epsi_1+\epsi_2}\\
e_{\epsi_3-\epsi_6}&e_{\epsi_3+\epsi_4}&-e_{\epsi_2+\epsi_3}&-e_{\epsi_1+\epsi_3}\\
e_{\epsi_4-\epsi_6}&0&-e_{\epsi_2+\epsi_4}&-e_{\epsi_1+\epsi_4}\\
\end{vmatrix}.
\end{equation*}

\lemmp{Let $\alpha\in\Phi^+\setminus\Bu$ and $\gamma\in\Phi^+$. If $\gamma\notin A(\alpha)$ then $[\Delta_{\alpha},e_{\gamma}]=0$\label{lemm:comm1_B_n_D_n}. If $\gamma\in A(\alpha)\setminus\Delta$ then $\alpha+\gamma\in\Phi^+$ and $[\Delta_{\alpha},e_{\gamma}]=\pm\Delta_{\alpha+\gamma}$.}{Let $\Phi=D_n$, $\alpha=\epsi_i-\epsi_j$, $i$ is even (all other cases can be considered similarly). One has
\begin{equation*}
\Delta_{\alpha}=\begin{vmatrix}
e_{\epsi_1-\epsi_j}&e_{\epsi_1+\epsi_i}&e_{\epsi_1+\epsi_{i-2}}&e_{\epsi_1+\epsi_{i-3}}&\ldots&
e_{\epsi_1+\epsi_3}&e_{\epsi_1+\epsi_2}&0\\
e_{\epsi_2-\epsi_j}&e_{\epsi_2+\epsi_i}&e_{\epsi_2+\epsi_{i-2}}&e_{\epsi_2+\epsi_{i-3}}&\ldots&
e_{\epsi_2+\epsi_3}&0&-e_{\epsi_1+\epsi_2}\\
e_{\epsi_3-\epsi_j}&e_{\epsi_3+\epsi_i}&e_{\epsi_3+\epsi_{i-2}}&e_{\epsi_3+\epsi_{i-3}}&\ldots&
0&-e_{\epsi_2+\epsi_3}&-e_{\epsi_1+\epsi_3}\\
e_{\epsi_4-\epsi_j}&e_{\epsi_4+\epsi_i}&e_{\epsi_4+\epsi_{i-2}}&e_{\epsi_4+\epsi_{i-3}}&\ldots&
-e_{\epsi_3+\epsi_4}&-e_{\epsi_2+\epsi_4}&-e_{\epsi_1+\epsi_4}\\
\vdots&\vdots&\vdots&\vdots&\ddots&\vdots&\vdots\\
e_{\epsi_{i-2}-\epsi_j}&e_{\epsi_{i-2}+\epsi_i}&0&-e_{\epsi_{i-3}+\epsi_{i-2}}&\ldots&
-e_{\epsi_3+\epsi_{i-2}}&-e_{\epsi_2+\epsi_{i-2}}&-e_{\epsi_1+\epsi_{i-2}}\\
e_{\epsi_{i-1}-\epsi_j}&e_{\epsi_{i-1}+\epsi_i}&-e_{\epsi_{i-2}+\epsi_{i-1}}&-e_{\epsi_{i-3}+\epsi_{i-1}}&\ldots&
-e_{\epsi_3+\epsi_{i-1}}&-e_{\epsi_2+\epsi_{i-1}}&-e_{\epsi_1+\epsi_{i-1}}\\
e_{\epsi_i-\epsi_j}&0&-e_{\epsi_{i-2}+\epsi_i}&-e_{\epsi_{i-3}+\epsi_i}&\ldots&
-e_{\epsi_3+\epsi_3}&-e_{\epsi_2+\epsi_i}&-e_{\epsi_1+\epsi_i}\\
\end{vmatrix}.
\end{equation*}

Suppose that $e_{\gamma}$, $\gamma\in\Phi^+$, does not commute with one of $e_{\delta}$'s involved in $\Delta_{\alpha}$. If $\gamma=\epsi_j+\epsi_k$ for some $k>j$ then, clearly, $[\Delta_{\alpha},e_{\gamma}]$ equals the minor obtained from $\Delta_{\alpha}$ by replacing all $e_{\epsi_s-\epsi_j}$ from the first column by $\epsi_s+\epsi_k$, but it is exactly $\Delta_{\epsi_i+\epsi_k}$. If $\gamma=\epsi_k+\epsi_j$ for some $i<k<j$ then $[\Delta_{\alpha},e_{\gamma}]$ equals the minor obtained from $\Delta_{\alpha}$ by replacing all $e_{\epsi_s-\epsi_j}$ from the first column by $-\epsi_s+\epsi_k$, but this minor is nothing but $-\Delta_{\epsi_i+\epsi_k}$. One can argue similarly for $\gamma=\epsi_{i-1}+\epsi_j$. In fact, we considered all possibilities for $\gamma\in A(\alpha)\setminus\Delta$.

Next, if $\gamma=\epsi_k+\epsi_j$ for some $k<i-1$, then $[\Delta_{\alpha},e_{\gamma}]$ equals the minor obtained from $\Delta_{\alpha}$ by replacing its first column by its $(i-k+1)$th column multiplied by $-1$, so $[\Delta_{\alpha},e_{\gamma}]=0$. If $\gamma=\epsi_k-\epsi_i$ for $k<i-1$ then $[\Delta_{\alpha},e_{\gamma}]=a+b$, where $a$ equals the minor obtained from $\Delta_{\alpha}$ by replacing its second column by its $(i-k+1)$th column multiplied by $(-1)$, and $b$ equals the minor obtained from $\Delta_{\alpha}$ by replacing its last row by its $k$th row multiplied by $(-1)$. Hence $a=b=0$. Arguing similarly, we conclude that $[\Delta_{\alpha},e_{\gamma}]=0$ for $\gamma=\epsi_k-\epsi_l$, $k<l<i$. For all other $\gamma\in\Phi^+$, $e_{\gamma}$ commutes with each $e_{\delta}$ involved in $\Delta_{\alpha}$. The proof is complete.}

\lemmp{Let $\alpha,~\gamma\in\Phi^+\setminus(\Bu\cap\Delta)$. If $\alpha+\gamma\notin\Bu$ then $[\Delta_{\alpha},\Delta_{\beta}]=0$. If $\alpha+\gamma=\epsi_i+\epsi_{i+1}\in\Bu$ then\label{lemm:comm2_B_n_D_n} there exists $a_{\alpha,\gamma}\in\Cp^{\times}$ such that
\begin{equation*}
[\Delta_{\alpha},\Delta_{\gamma}]=\begin{cases}a_{\alpha,\gamma}P_{\epsi_1+\epsi_2}^3&\text{if }i=1,\\
a_{\alpha,\gamma}P_{\epsi_i+\epsi_{i+1}}^3P_{\epsi_{i-2}+\epsi_{i-1}}&\text{otherwise}.\\
\end{cases}
\end{equation*}\vspace{-0.1cm}}
{Suppose $\Phi=D_n$ (the proof for $B_n$ is similar). Pick two roots $\alpha$, $\gamma\in\Phi^+\setminus(\Bu\cap\Delta)$. First, we claim that $[\Delta_{\alpha},\Delta_{\gamma}]=0$ if $\alpha+\gamma\notin\Bu$. We will consider the case $\alpha=\epsi_i-\epsi_j$ for even $i$ (all other cases can be considered similarly). Assume that there exists $e_{\delta}$ involved in $\Delta_{\gamma}$, which do not commute with $\Delta_{\alpha}$. It follows from the definitions of $A(\alpha)$ and $\Delta_{\gamma}$ that either $\row(\gamma)>i$ and $\col(\gamma)\neq-j$, or $\gamma\in\{\epsi_{i-1}+\epsi_j\}\cup\bigcup_{k=i+1}^{j-1}\{\epsi_k+\epsi_j\}$. If there exists $e_{\delta}$ involved in $\Delta_{\alpha}$, which do not commute with $\Delta_{\gamma}$, then, by the definitions of $A(\gamma)$ and $\Delta_{\alpha}$, one has $\gamma=\epsi_{i-1}+\epsi_j$. In this case, denote by $\Delta^k$, $1\leq k\leq i$, the $(i-1)\times(i-1)$ minor obtained from $\Delta_{\alpha}$ (or, equivalently, from $\Delta_{\gamma}$) by deleting its first row and its $k$th column. One has
\begin{equation*}
\begin{split}
[\Delta_{\alpha},\Delta_{\gamma}]&=\left[\sum\nolimits_{k=1}^i(-1)^ke_{\epsi_k-\epsi_j}\Delta^k,
\sum\nolimits_{l=1}^i(-1)^le_{\epsi_l+\epsi_j}\Delta^l\right]\\
&=\sum_{1\leq k<l\leq i}(-1)^{k+l}[e_{\epsi_k-\epsi_j},e_{\epsi_l+\epsi_j}]\Delta^k\Delta^l+
\sum_{1\leq l<k\leq i}(-1)^{k+l}[e_{\epsi_k-\epsi_j},e_{\epsi_l+\epsi_j}]\Delta^k\Delta^l\\
&=\sum_{1\leq k<l\leq i}(-1)^{k+l}[e_{\epsi_k-\epsi_j},e_{\epsi_l+\epsi_j}]\Delta^k\Delta^l-
\sum_{1\leq k<l\leq i}(-1)^{k+l}[e_{\epsi_k+\epsi_j},e_{\epsi_l-\epsi_j}]\Delta^k\Delta^l\\
&=-\sum_{1\leq k<l\leq i}(-1)^{k+l}e_{\epsi_k+\epsi_l}\Delta^k\Delta^l
+\sum_{1\leq k<l\leq i}(-1)^{k+l}e_{\epsi_k+\epsi_l}\Delta^k\Delta^l=0.
\end{split}
\end{equation*}
Thus, if $\alpha+\gamma\notin\Bu$ then $[\Delta_{\alpha},\Delta_{\gamma}]=0$, as required.

Now assume that $\alpha+\gamma=\epsi_i+\epsi_{i+1}\in\Bu$. Then $[\Delta_{\alpha},\Delta_{\gamma}]\in Z(\nt)$. Indeed, given $\delta\in\Phi^+$, we have
$$[[\Delta_{\alpha},\Delta_{\gamma}],e_{\delta}]=[\Delta_{\alpha},[\Delta_{\gamma},e_{\delta}]]+
[[\Delta_{\alpha},e_{\delta}],\Delta_{\gamma}].$$
According to Lemma~\ref{lemm:comm1_B_n_D_n}, the RHS is a linear combination of $[\Delta_{\alpha+\delta},\Delta_{\gamma}]$ and $[\Delta_{\gamma+\delta},\Delta_{\alpha}]$. (For $\mu\notin\Phi^+$, we put by definition $\Delta_{\mu}=0$.) But if $\alpha+\gamma\in\Bu$ and $\delta\in\Phi^+$ then $\alpha+\gamma+\delta\notin\Bu$, so both $[\Delta_{\alpha+\delta},\Delta_{\gamma}]$ and $[\Delta_{\gamma+\delta},\Delta_{\alpha}]$ are zero by the above. Thus, $[\Delta_{\alpha},\Delta_{\gamma}]\in Z(\nt)$.

Recall that the Cartan subalgebra $\htt$ of $\gt$ acts on $U(\nt)$. Given $\delta\in(\Ro_k\cup\Ro_{k-1})\setminus\Bu$, $k$ even, it is clear from the definition that $\Delta_{\delta}$ is an $\htt$-weight element of weight $$\vpi(\delta)=2\epsi_1+\ldots+2\epsi_{k-2}+\epsi_{k-1}+\epsi_k+\delta.$$ Hence, $[\Delta_{\alpha},\Delta_{\gamma}]$ is an $\htt$-weight element of $Z(\nt)$ of weight $$\vpi(\alpha)+\vpi(\gamma)=2\epsi_1+\ldots+2\epsi_{i-2}+3\epsi_{i-1}+3\epsi_i.$$ It follows from (\ref{formula:weights_of_xi_i}) that $P_{\epsi_i+\epsi_{i+1}}^3P_{\epsi_{i-2}+\epsi_{i-1}}$ ($P_{\epsi_1+\epsi_2}$ for $i=1$) has the same weight. Since $Z(\nt)$ is a direct sum of $1$-dimensional $\htt$-weight subspaces \cite[Theorem 6]{Kostant2}, \cite[Lemma 4.4]{Joseph1}, it remains to check that $[\Delta_{\alpha},\Delta_{\gamma}]$ is nonzero.

To do this, denote by $\Delta_{\alpha}^k$ (respectively, by $\Delta_{\gamma}^k$), $1\leq k\leq i+1$, the minor  obtained from $\Delta_{\alpha}$ (respectively, from $\Delta_{\gamma}$) by deleting its first row and its $k$th column. We may assume without loss of generality that $\alpha=\epsi_i-\epsi_j$ and $\gamma=\epsi_{i+1}+\epsi_j$ for some $j\geq i+1$. Then
\begin{equation*}
\begin{split}
[\Delta_{\alpha},\Delta_{\gamma}]&=\left[\sum_{k=1}^{i+1}(-1)^ke_{\epsi_k-\epsi_j}\Delta_{\alpha}^k,
\sum_{l=1}^{i+1}(-1)^le_{\epsi_l+\epsi_j}\Delta_{\gamma}^l\right]\\
&=\sum_{k\neq l}(-1)^{k+l}[e_{\epsi_k-\epsi_j},e_{\epsi_l+\epsi_j}]\Delta_{\alpha}^k\Delta_{\gamma}^l\\
&=\pm e_{\epsi_i+\epsi_{i+1}}\Delta_{\alpha}^{i+1}\Delta_{\gamma}^i+\sum_{k\neq l}(-1)^{k+l}[e_{\epsi_k-\epsi_j},e_{\epsi_l+\epsi_j}]\Delta_{\alpha}^k\Delta_{\gamma}^l.
\end{split}
\end{equation*}
(The first summand corresponds to $k=i+1$, $l=i$, and the second one corresponds to all other $k\neq l$.) The first summand contains exactly one term of the form $\pm e_{\epsi_i+\epsi_{i+1}}^3P_{\epsi_{i-2}+\epsi_{i-1}}^2$ (if $i=1$, then exactly one term of the form $\pm e_{\epsi_i+\epsi_{i+1}}^3$). On the other hand, $e_{\epsi_i+\epsi_{i+1}}$ enters each term of the second summand with degree at most $2$. We conclude that $[\Delta_{\alpha},\Delta_{\gamma}]\neq0$. The proof is complete.}

\corop{Given\label{coro:comm2_B_n_D_n} $\alpha\in\Phi^+\setminus\Bu$\textup, there exists the unique nonzero element $P_{\alpha}$ of $U(\nt)$ such that $\Delta_{\alpha}=P_{\alpha}P_{\beta}$\textup, where $\beta=\epsi_i+\epsi_{i+1}$ if $\row(\alpha)=i$ is odd or $\row(\alpha)=i+1$ is even. Further\textup, let $\alpha,~\gamma\in\Phi^+\setminus\Bu$. If $\alpha+\gamma\notin\Bu$ then $[P_{\alpha},P_{\gamma}]=0$. If $\alpha+\gamma=\epsi_i+\epsi_{i+1}\in\Bu$ then
\begin{equation*}
[P_{\alpha},P_{\gamma}]=\begin{cases}a_{\alpha,\gamma}P_{\epsi_1+\epsi_2}&\text{if }i=1,\\
a_{\alpha,\gamma}P_{\epsi_i+\epsi_{i+1}}P_{\epsi_{i-2}+\epsi_{i-1}}&\text{otherwise}.\\
\end{cases}
\end{equation*}
Here $a_{\alpha,\gamma}$ are as in Lemma~\textup{\ref{lemm:comm2_B_n_D_n}}.
}{We say that a square matrix is \emph{skew-antisymmetric} if it is skew-symmetric with respect to the antidiagonal. Note that the Pfaffian of a skew-antisymmetric matrix is well defined, and that the matrix $\Uu$ is skew-antisymmetric. Given $\alpha\in\Phi^+\setminus\Bu$, denote by $\Uu_{\alpha}$ the smallest skew-antisymmetric submatrix of $\Uu$ containing the rows $R(\alpha)$ and the columns $C(\alpha)$. Let $P_{\alpha}'$ be the Pfaffian of $\Uu_{\alpha}$. It follows from Cayley's results \cite{Cayley1} (see also \cite[(2.24)]{Heymans1}) that $\Delta_{\alpha}=\pm P_{\alpha}'P_{\beta}$, where $\beta=\epsi_i+\epsi_{i+1}$ if $\row(\alpha)=i$ is odd or $\row(\alpha)=i+1$ is even. Putting $P_{\alpha}=\pm P_{\alpha}'$, we see that $\Delta_{\alpha}=P_{\alpha}P_{\beta}$. All other claims of the corollary follow immediately from Lemma~\ref{lemm:comm2_B_n_D_n} and the fact that the algebra $U(\nt)$ does not contain zero divisors.}

We are ready to prove the second result of this subsection.

\theop{Suppose\label{theo:Weyl_quotient} $\Phi$ is of type $B_n$ or $D_n$. Let $J$ be a centrally generated primitive ideal of $U(\nt)$\textup, and let $\wh x$ denote the image in $U(\nt)/J$ of an element $x\in U(\nt)$ under the canonical projection. Then\textup, for every odd $i$ and every $\alpha\in\Ro_i$\textup, there exists a constant $a_{\alpha}$ such that the elements $p_{\alpha}=\wh P_{\alpha}$ and $q_{\alpha}=a_{\alpha}\wh P_{\epsi_i+\epsi_{i+1}-\alpha}$ satisfy the following conditions\textup:
\begin{equation*}
\begin{split}\text{\textup{i) }}&\text{all $p_{\alpha}$ and $q_{\alpha}$ generate the quotient $U(\nt)/J$ as an algebra\textup;}\\
\text{\textup{ii) }}&[p_{\alpha},q_{\alpha}]=1\text{ for all }\alpha,~[p_{\alpha},p_{\beta}]=[q_{\alpha},q_{\beta}]=0\text{ for all }\alpha,~\beta,\text{ and }[p_{\alpha},q_{\beta}]=0\text{ for }\alpha\neq\beta.\\
\end{split}
\end{equation*}}{Consider the case $\Phi=D_n$ (the proof for $B_n$ is similar). Put $\Au=U(\nt)/J$. Clearly, $\Au$ is generated as an algebra by $\wh e_{\alpha}$, $\alpha\in\Phi^+$.

First, we claim that $\Au$ is generated by $\wh e_{\alpha}$, $\alpha\in\Phi^+\setminus(\Bu\cap\Delta)$. Indeed, recall the notion of $\wt\nt$, $\wt\Phi$, $J_0$, $\wt\Bu$, $\wt P_{\beta}$ for $\beta\in\wt\Bu$, and $\wt J_c$ from the proof of Proposition~\ref{prop:J_c_primitive_Bn_Dn}. In our case, $J=J_c$, and we denote for simplicity $\wt J=\wt J_c$. The epimorphism $r\colon U(\nt)\sur U(\wt\nt_e)\otimes\Au_s$ induces the isomorphism of algebras $$\wh r\colon\Au=U(\nt)/J\to r(U(\nt))/r(J),~\wh r(\wh x)=r(x)+r(J)\text{ for }x\in U(\nt).$$ It follows from the proof of Proposition~\ref{prop:J_c_primitive_Bn_Dn} that $r(U(\nt))/r(J)\cong(U(\nt)/J_0)/r(J)$ is isomorphic to $(U(\wt\nt)/\wt J)\otimes\Au_s$. Slightly abusing the notation, we denote the image of an element $x\in U(\wt\nt)$ in the quotient algebra $U(\wt\nt)/\wt J$ again by $\wh x$. Using the induction on $\rk\Phi$, we may assume that $U(\wt\nt)/\wt J$ is generated by $\wh e_{\alpha}$ for $\alpha\in\wt\Phi^+\setminus(\wt\Bu\cap\Delta)$ (the base $\rk\Phi\leq3$ can be easily checked).

Hence, $\Au$ is generated as an algebra by $(\wh r)^{-1}(\wh e_{\alpha}\otimes1)$ for $\alpha\in\wt\Phi^+\setminus(\wt\Bu\cap\Delta)$, and by $\wh e_{\alpha}=(\wh r)^{-1}(1\otimes\overline{e}_{\alpha})$ for $\alpha\in(\Ro_1\cup\Ro_2)\setminus\{\beta_1,~\beta_2\}$. As we mentioned in the proof of Proposition~\ref{prop:J_c_primitive_Bn_Dn}, for $\alpha\in\wt\Phi$,
$$r(e_{\alpha})=e_{\alpha}\otimes1+1\otimes(\text{linear combination of $\overline{e}_{\alpha+\gamma}\overline{e}_{\beta_1-\gamma},~\gamma\in(\Ro_1\cup\Ro_2)\setminus\{\beta_1,~\beta_2\}$}),$$ and, consequently,
$$\wh r(\wh e_{\alpha})=\wh e_{\alpha}\otimes1+1\otimes(\text{linear combination of $\overline{e}_{\alpha+\gamma}\overline{e}_{\beta_1-\gamma},~\gamma\in(\Ro_1\cup\Ro_2)\setminus\{\beta_1,~\beta_2\}$}).$$
But both $\alpha+\gamma$ and $\beta_1-\gamma$ belong to $(\Ro_1\cup\Ro_2)\setminus\{\beta_1,~\beta_2\}$, so
$$(\wh r)^{-1}(\wh e_{\alpha}\otimes1)=\wh e_{\alpha}-(\text{linear combination of $\wh e_{\alpha+\gamma}\wh e_{\beta_1-\gamma},~\gamma\in(\Ro_1\cup\Ro_2)\setminus\{\beta_1,~\beta_2\}$}).$$
We conclude that $\Au$ is generated by $\wh e_{\alpha}$ for $\alpha\in(\wt\Phi^+\setminus(\wt\Bu\cap\Delta))\cup((\Ro_1\cup\Ro_2)\setminus\{\beta_1,~\beta_2\})\subset\Phi^+\setminus(\Bu\cap\Delta)$, as required.

Second, assume that $\beta=\epsi_i+\epsi_{i+1}\in\Bu\setminus\Delta$ for odd $i>1$. (Clearly, $\wh e_{\epsi_1+\epsi_2}=c_{\epsi_1+\epsi_2}$ in $\Au$.) Recall that there is a natural partial order on $\Phi^+$: $\alpha>\gamma$ if $\alpha-\gamma$ is a sum of positive roots. We have
$$P_{\beta}=e_{\beta}P_{\epsi_{i-2}+\epsi_{i-1}}+\text{terms containing only $e_{\alpha}$ for }\alpha>\beta.$$
It follows that in $\Au$ we can write $\wh e_{\beta}$ as a polynomial in $\wh e_{\alpha}$ for $\alpha\in\Phi^+\setminus\Bu$. In other words, $\Au$~is generated as an algebra by $\wh e_{\alpha}$ for $\alpha\in\Phi^+\setminus\Bu$.

Finally, let $\alpha\in\Phi^+\setminus\Bu$. Assume that $\row(\alpha)=i$ is odd or $\row(\alpha)=i+1$ is even. If $i>1$ then $\Delta_{\alpha}=\pm e_{\alpha}e_{\epsi_1+\epsi_2}$. It follows again from \cite{Cayley1} (see also \cite[2.23]{Heymans1}) that, for $i>1$,
\begin{equation*}
\Delta_{\alpha}=\pm e_{\alpha}P_{\epsi_i+\epsi_{i+1}}P_{\epsi_{i-2}+\epsi_{i-1}}+\text{terms containing only $e_{\gamma}$ with $\gamma>\alpha$.}
\end{equation*}
This implies that in $\Au$ one can write $\wh e_{\alpha}$ as a polynomial in $\wh\Delta_{\gamma}$ for $\gamma\in\Phi^+\setminus\Bu$. Since $\wh P_{\beta}=c_{\beta}\neq0$ in $\Au$ for each $\beta\in\Bu\setminus\Delta$, Corollary~\ref{coro:comm2_B_n_D_n} shows that once can write $\wh e_{\alpha}$ as a polynomial in $\wh P_{\gamma}$ with $\gamma\in\Phi^+\setminus\Bu$. Thus, $\wh\Delta_{\gamma}$ with $\gamma\in\Phi^+\setminus\Bu$ generate $\Au$ as an algebra. Now, given $\alpha\in\Phi^+\setminus\Bu$ with odd $\row(\alpha)=i$, put $a_{\alpha}=a_{\alpha,\epsi_i+\epsi_{i+1}-\alpha}^{-1}$. To conclude the proof, it remains to apply Corollary~\ref{coro:comm2_B_n_D_n} again.}

\nota{The fact that $U(\nt)/J$ is generated as an algebra by the elements $\wh e_{\alpha}$ for $\alpha\in\Phi^+\setminus\Bu$ plays a crucial role in the proof of Corollary~\ref{coro:comm2_B_n_D_n}. For $A_{n-1}$ and $C_n$ the proof of this fact was more or less obvious, but for $B_n$ and $D_n$ this is non-trivial, because the explicit form of $P_{\beta}$ for $\beta\in\Bu\cap\Delta$ is rather complicated. That's why we involve the epimorphism $r$ into the picture and use the induction on $\rk\Phi$. In fact, a posteriori one can easily give a proof of Proposition~\ref{prop:J_c_primitive_Bn_Dn} based on Corollary~\ref{coro:comm2_B_n_D_n} in a similar way as it was done for $A_{n-1}$ and $C_n$ in \cite{IgnatyevPenkov1}.}

\sect{Infinite-dimensional\label{sect:infinite_dim} case}

\sst{The center of $U(\nt)$} Let \label{sst:centers_ifd} $\slt_{\infty}(\Cp)$, $\ot_{\infty}(\Cp)$, $\spt_{\infty}(\Cp)$ be the three simple complex finitary countable dimensional Lie algebras as classified by A.~Baranov \cite{Baranov1}. Each of them can be described as follows (see for example \cite{DimitrovPenkov2}). Consider an infinite chain of inclusions $$\gt_1\subset\gt_2\subset\ldots\subset\gt_n\subset\ldots$$ of simple Lie algebras, where $\rk{\gt_n}=n$ and all $\gt_n$ are of the same type $A$, $B$, $C$ or $D$. Then the union $\gt=\bigcup\gt_n$ is isomorphic to $\slt_{\infty}(\Cp)$, $\ot_{\infty}(\Cp)$ or $\spt_{\infty}(\Cp)$. It is always possible to choose nested Cartan subalgebras $\htt_n\subset\gt_n$, $\htt_n\subset\htt_{n+1}$, so that each root space of $\gt_n$ is mapped to exactly one root space of $\gt_{n+1}$. The union $\htt=\bigcup\htt_n$ acts semisimply on $\gt$, and is by definition a \emph{splitting Cartan subalgebra} of~$\gt$. We have a root decomposition $$\gt=\htt\oplus\bigoplus_{\alpha\in\Phi}\gt^{\alpha}$$ where $\Phi$ is \emph{the root system of} $\gt$ with respect to $\htt$ and $\gt^{\alpha}$ are the \emph{root spaces}. The root system $\Phi$ is simply the union of the root systems of $\gt_n$ and equals one of the following infinite root systems:
\begin{equation*}
\begin{split}
A_{\infty}&=\pm\{\epsi_i-\epsi_j,~i,j\in\Zp_{>0},~i<j\},\\
B_{\infty}&=\pm\{\epsi_i-\epsi_j,~i,j\in\Zp_{>0},~i<j\}\cup\pm\{\epsi_i+\epsi_j,~i,j\in\Zp_{>0},~i<j\}\cup\pm\{\epsi_i,~i\in\Zp_{>0}\},\\
C_{\infty}&=\pm\{\epsi_i-\epsi_j,~i,j\in\Zp_{>0},~i<j\}\cup\pm\{\epsi_i+\epsi_j,~i,j\in\Zp_{>0},~i<j\}\cup\pm\{2\epsi_i,~i\in\Zp_{>0}\},\\
D_{\infty}&=\pm\{\epsi_i-\epsi_j,~i,j\in\Zp_{>0},~i<j\}\cup\pm\{\epsi_i+\epsi_j,~i,j\in\Zp_{>0},~i<j\}.
\end{split}
\end{equation*}

A \emph{splitting Borel subalgebra} of $\gt$ is a subalgebra $\bt$ such that for every~$n$, $\bt_n=\bt\cap\gt_n$ is a Borel subalgebra of $\gt_n$. It is well-known that any splitting Borel subalgebra is conjugate via $\mathrm{Aut}\,\gt$ to a splitting Borel subalgebra containing $\htt$. Therefore, in what follows we restrict ourselves to considering only such Borel subalgebras $\bt$.

Recall \cite{DimitrovPenkov1} that a linear order on $\{0\}\cup\{\pm\epsi_i\}$ is $\Zp_2$-\emph{linear} if multiplication by $-1$ reverses the order. By \cite[Proposition 3]{DimitrovPenkov1}, there exists a bijection between splitting Borel subalgebras of $\gt$ containing $\htt$ and certain linearly ordered sets as follows.
\begin{equation*}\predisplaypenalty=0
\begin{split}
&\text{For }A_{\infty}\text{: linear orders on }\{\epsi_i\};\\
&\text{for }B_{\infty}\text{ and }C_{\infty}\text{: }\Zp_2\text{-linear orders on }\{0\}\cup\{\pm\epsi_i\};\\
&\text{for }D_{\infty}\text{: }\Zp_2\text{-linear orders on }\{0\}\cup\{\pm\epsi_i\}\text{ with the property that}\\
&\text{a minimal positive element (if it exists) belongs to }\Zp_{>0}.
\end{split}
\end{equation*}
In the sequel we denote these linear orders by $\prec$. To write down the above bijection, denote $\teta_i=\epsi_i$, if $i\succ0$, and $\teta_i=-\epsi_i$, if $\epsi_i\prec0$ (for $A_{\infty}$, $\teta_i=\epsi_i$ for all $i$). Then put $\bt=\htt\oplus\nt$, where $\nt=\bigoplus\limits_{\alpha\in\Phi^+}\gt^{\alpha}$ and
\begin{equation*}\predisplaypenalty=0
\begin{split}
A_{\infty}^+&=\{\teta_i-\teta_j,~i,j\in\Zp_{>0},~\teta_i\succ\teta_j\},\\
B_{\infty}^+&=\{\teta_i-\teta_j,~i,j\in\Zp_{>0},~\teta_i\succ\teta_j\}\cup\{\teta_i+\teta_j,~i,j\in\Zp_{>0},
~\teta_i\succ\teta_j\}\cup\{\teta_i,~i\in\Zp_{>0}\},\\
C_{\infty}^+&=\{\teta_i-\teta_j,~i,j\in\Zp_{>0},~\teta_i\succ\teta_j\}\cup\{\teta_i+\teta_j,~i,j\in\Zp_{>0},~
\teta_i\succ\teta_j\}\cup\{2\teta_i,~i\in\Zp_{>0}\},\\
D_{\infty}^+&=\{\teta_i-\teta_j,~i,j\in\Zp_{>0},~\teta_i\succ\teta_j\}\cup\{\teta_i+\teta_j,~i,j\in\Zp_{>0},~\teta_i\succ\teta_j\}.\\
\end{split}
\end{equation*}

Our goal in this subsection is to recall the description of the center $Z(\nt)$ of the enveloping algebra $U(\nt)$ from \cite{IgnatyevPenkov1}. Fix $\nt$, i.e., fix an order $\prec$ as above. Define the subset $\Nu\subseteq\Zp_{>0}$ by setting $\Nu=\bigcup_{k\geq0}\Nu_k$, where $\Nu_0=\varnothing$ and $\Nu_k$ for $k\geq1$ is defined inductively in the following table.
\begin{center}
\begin{tabular}{|l|l|}
\hline
$\Phi$&$\Nu_k$\\
\hline\hline
$A_{\infty}$&$\Nu_{k-1}\cup\{i,j\}$ if there exists a maximal element $\teta_i$\\
&and a minimal element $\teta_j$ of $\{\teta_s,~s\in\Zp_{>0}\setminus\Nu_{k-1}\}$,\\
&$\Nu_{k-1}$ otherwise\\
\hline
$C_{\infty}$&$\Nu_{k-1}\cup\{i\}$ if there exists a maximal element $\teta_i$ of $\{\teta_s,~s\in\Zp_{>0}\setminus\Nu_{k-1}\}$,\\
&$\Nu_{k-1}$ otherwise\\
\hline
$B_{\infty}$,&$\Nu_{k-1}\cup\{i,j\}$ if there exists a maximal element $\teta_i$ of $\{\teta_s,~s\in\Zp_{>0}\setminus\Nu_{k-1}\}$\\
$D_{\infty}$&and a maximal element $\teta_j$ of $\{\teta_s,~s\in\Zp_{>0}\setminus\left(\Nu_{k-1}\cup\{i\}\right)\}$,\\
&$\Nu_{k-1}$ otherwise\\
\hline
\end{tabular}
\end{center}

\exam{i) Let $\Phi=A_{\infty}$. If $\epsi_1\succ\epsi_3\succ\ldots\succ\epsi_4\succ\epsi_2$, then $\Nu=\Zp_{>0}$. If $\epsi_1\succ\epsi_2\succ\epsi_3\succ\ldots$, then $\Nu=\varnothing$. ii) Let $\Phi\neq A_{\infty}$ and $\epsi_1\succ\epsi_2\succ\ldots\succ0\succ\ldots\succ-\epsi_2\succ-\epsi_1$. Then $\Nu=\Zp_{>0}$.}

Now we can define the (possibly infinite) Kostant cascade corresponding to $\nt$. Namely, to each $\Nu_k$ such that $\Nu_{k-1}\subsetneq\Nu_k$, we assign the root
\begin{equation*}
\beta_k=\begin{cases}\teta_i-\teta_j,&\text{if $\Phi=A_{\infty}$ and $\Nu_k\setminus\Nu_{k-1}=\{i,j\}$, $i\succ j$,}\\
\teta_i+\teta_j,&\text{if $\Phi=B_{\infty}$ or $D_{\infty}$ and $\Nu_k\setminus\Nu_{k-1}=\{i,j\}$, $i\succ j$,}\\
2\teta_i,&\text{if $\Phi=C_{\infty}$ and $\Nu_k\setminus\Nu_{k-1}=\{i\}$,}\\
\end{cases}
\end{equation*}
and put $$\Bu=\{\beta_k,~k\geq1,~\Nu_{k-1}\subsetneq\Nu_k\}.$$
Note that $\Bu$ is a strongly orthogonal subset of $\Phi^+$; however it is not necessarily maximal with this property.

\defi{The subset $\Bu$ is called the \emph{Kostant cascade} corresponding to $\nt$.}
\exam{i) If $\Phi=A_{\infty}$ and $\epsi_1\succ\epsi_3\succ\ldots\succ\epsi_4\succ\epsi_2$, then $$\Bu=\{\epsi_1-\epsi_2,~\epsi_3-\epsi_4,~\epsi_5-\epsi_6,~\ldots\}.$$ ii) If $\Phi\neq A_{\infty}$ and $\epsi_1\succ\epsi_2\succ\ldots\succ0\succ\ldots\succ-\epsi_2\succ-\epsi_1$, then
\begin{equation*}
\Bu=\begin{cases}\{\epsi_1+\epsi_2,~\epsi_3+\epsi_4,~\epsi_5+\epsi_6,~\ldots\}&\text{for $B_{\infty}$ and $D_{\infty}$},\\
\{2\epsi_1,~2\epsi_2,~2\epsi_3,~\ldots\}&\text{ for $C_{\infty}$}.
\end{cases}
\end{equation*}}

To each finite non-empty subset $M\subset\Zp_{>0}$, one can assign a root subsystem $\Phi_M$ of $\Phi$ and a subalgebra~$\nt_M$ of $\nt$ by putting
\begin{equation*}
\begin{split}
\Phi_M&=\Phi\cap\langle\epsi_i,~i\in M\rangle_{\Rp},\\
\nt_M&=\bigoplus_{\alpha\in\Phi_M^+}\gt^{\alpha},~\Phi_M^+=\Phi_M\cap\Phi^+.
\end{split}
\end{equation*}
Then the subsystem $\Phi_M$ is isomorphic to the root system $\Phi_n$ of $\gt_n$ for $n=|M|$; we denote this isomorphism by $j_M\colon\Phi_n\to\Phi_M,~\epsi_i\mapsto\teta_{a_i}$, where $M=\{a_1,\ldots,a_n\}$, $\teta_{a_1}\succ\ldots\succ\teta_{a_n}$. Besides, $\nt_M$ is isomorphic as a Lie algebra to the maximal nilpotent subalgebra $\nt_n$ of $\gt_n$ considered in the previous section. Note also that $\nt=\ilm\nt_M$. Here, for $M\subseteq M'$, the monomorphism ${i_{M,M'}\colon\nt_M\hookrightarrow\nt_{M'}}$ is just the inclusion. Further, it is easy to see that there exist isomorphisms $\phi_M\colon\nt_n\to\nt_M,~M\subset\Zp_{>0},~n=|M|$, such that, for $M\subseteq M'$, $i_{M,M'}\circ\phi_M$ is just the restriction of $\phi_{M'}$ to $\nt_n\subset\nt_{n'}$, $n'=|M'|$, and, for $\alpha\in\Phi_n^+$, $\phi_M(e_{\alpha})$ is a root vector corresponding to the root $j_M(\alpha)$; we denote it by $f_{j_M(\alpha)}$.

We are now ready to write down a set of generators of $Z(\nt)$. Namely, suppose that $\beta=\beta_k\in\Bu$ for some $k\geq1$. Let $M$ be a finite subset of $\Zp_{>0}$ such that $\Nu_k\subseteq M$. The isomorphism $\phi_M$ gives rise to an isomorphism $U(\nt_n)\to U(\nt_M)$, $n=|M|$. We denote the respective images of $\Delta_{j_M^{-1}(\beta)}$ and $P_{j_M^{-1}(\beta)}$ (as elements of $U(\nt_n)$ whenever defined) in $U(\nt_M)$ by $\Delta_{\beta}$ and $P_{\beta}$. Then for $A_{\infty}$ (respectively, for~$C_{\infty}$), $\Delta_{\beta}\in U(\nt_M)$ is given by formula (\ref{formula:Delta_i_A_n}) (respectively, (\ref{formula:Delta_i_C_n})) for $j_M^{-1}(\beta)=\epsi_i-\epsi_{n-i+1}$ (respectively, for $j_M^{-1}(\beta)=2\epsi_i$) with $f_{j_M(\alpha)}$ instead of $e_{\alpha}$. Similarly, for $B_{\infty}$ and $D_{\infty}$, $P_{\beta}\in U(\nt_M)$ is given by formula (\ref{formula:Delta_i_D_n_pf}) for $j_M^{-1}(\beta)=\epsi_i+\epsi_{i+1}$ with $f_{j_M(\alpha)}$ instead of~$e_{\alpha}$. It is important that $\Delta_{\beta}$, $P_{\beta}\in U(\nt_M)$ depend only on $\beta$ but not on $M$. Moreover, it is clear from the finite-dimensional theory that $\Delta_{\beta}$ (respectively, $P_{\beta}$) belong to the center of $U(\nt)$ for $A_{\infty}$ and $C_{\infty}$ (respectively, for $B_{\infty}$ and $D_{\infty}$).

The following theorem was proved in \cite[Theorem 2.6]{IgnatyevPenkov1}.

\mtheo{If \label{theo:center_ifd} $\Phi=A_{\infty}$\textup{,} $C_{\infty}$ \textup{(}respectively\textup{,} $\Phi=B_{\infty}$\textup{,} $D_{\infty}$\textup{)}\textup{,} then $\Delta_{\beta}$ \textup{(}res\-pec\-tively\textup{,}~$P_{\beta}$\textup{)}\textup{,} $\beta\in\Bu$\textup{,} generate $Z(\nt)$ as an algebra. In particular\textup{,} $Z(\nt)$ is a polynomial ring in $|\Bu|$ variables.}

\sst{Centrally generated ideals of $U(\nt)$} Throughout\label{sst:ideals_ifd} this subsection we use the notation from Sub\-sec\-tion~\ref{sst:centers_ifd}. We now restrict ourselves to the case $\Nu=\Zp_{>0}$. This means that, up to isomorphism, $\nt$ can be chosen to correspond to the linear order $\epsi_1\succ\epsi_3\succ\epsi_5\succ\ldots\succ\epsi_6\succ\epsi_4\succ\epsi_2$ for $A_{\infty}$ (respectively, to the linear order
$\epsi_1\succ\epsi_2\succ\epsi_3\succ\ldots\succ0\succ\ldots\succ-\epsi_2\succ-\epsi_1$ for all other root systems). In particular, $\teta_i=e_i$ for all $i\in\Zp_{>0}$, and $f_{\alpha}=e_{\alpha}$ for all $M\subset\Zp_{>0}$, $\alpha\in\Phi_M^+$. For $\alpha\in\Phi^+$, denote by $e_{\alpha}^*$ the linear form on $\nt$ such that $e_{\alpha}^*(e_{\beta})=\delta_{\alpha,\beta}$ (the Kronecker delta) for all $\beta\in\Phi^+$. In this subsection we describe all centrally generated primitive ideals of $U(\nt)$ for $B_{\infty}$ and $D_{\infty}$.

For our choice of $\nt$, the Kostant cascade has the following form:
\begin{equation*}
\Bu=\begin{cases}
\{\epsi_1-\epsi_2,~\epsi_3-\epsi_4,~\ldots\}&\text{for }A_{\infty},\\[-2pt]
\{\epsi_1+\epsi_2,~\epsi_3+\epsi_4,~\ldots\}&\text{for }B_{\infty}\text{ and }D_{\infty},\\[-2pt]
\{2\epsi_1,~2\epsi_2,~\ldots\}&\text{for }C_{\infty}.\\
\end{cases}
\end{equation*}
The forms $$f_{\xi}=\sum_{\beta\in\Bu}\xi(\beta)e_{\beta}^*\in\nt^*$$ for all maps $\xi\colon\Bu\to\Cp^{\times}$ are by definition the \emph{Kostant forms} on $\nt$.

Our goal is to construct for $B_{\infty}$ and $D_{\infty}$ a partial Dixmier map, which attaches to each Kostant form a primitive ideal of $U(\nt)$ (for $A_{\infty}$ and $C_{\infty}$ it was done in \cite{IgnatyevPenkov1}). As in Subsection~\ref{sst:ideals_fd}, define the maps $\row\colon\Phi\to\Zp$ and $\col\colon\Phi\to\Zp$ by putting
\begin{equation*}
\begin{split}
&\row(\epsi_i-\epsi_j)=\row(\epsi_i+\epsi_j)=\row(2\epsi_i)=\row(\epsi_i)=i,\\
&\col(\epsi_i+\epsi_j)=\col(2\epsi_j)=-j,\\
\end{split}
\end{equation*}
and set $\Ro_k=\{\alpha\in\Phi^+\mid\row(\alpha)=k\}$. Put $\pt=\langle e_{\alpha},~\alpha\in\Phi^+\setminus\Mo\rangle_{\Cp}$, where
\begin{equation*}
\Mo=\begin{cases}
\{\epsi_i-\epsi_j,~\text{$i$ odd, $j$ even, $j<i$}\}&\text{for }A_{\infty},\\[-2pt]
\Ro_i,~\text{$i$ even}&\text{for $B_{\infty}$ and $D_{\infty}$},\\[-2pt]
\{\epsi_i-\epsi_j,~1\leq i<j\leq n\}&\text{for }C_{\infty}.\\
\end{cases}
\end{equation*}
Put also $\pt_n=\pt\cap\nt_n$, where $\nt_n=\nt_M$ for $M=\{1,\ldots,n\}$. Fix a Kostant form $f=f_{\xi}$. By\break \cite[Theo\-rem~1.1]{Ignatyev1}, $\pt_n$ is a polarization of $\nt_n$ at the linear form $f_n=f\mathbin{\mid}_{\nt_n}$. Thus, $\pt=\ilm\pt_n$ is a~polarization of $\nt$ at~$f$. Moreover, denote
\begin{equation}
V_{\xi}=U(\nt)\otimes_{U(\pt)}W,~V_{\xi}^n=U(\nt_n)\label{formula:V_V_m_ifd}\otimes_{U(\pt_n)}W^n,
\end{equation}
where $W$ (respectively, $W^n$) is the one-dimensional representation of $\pt$ (respectively, of $\pt_n$) given by $x\mapsto f_{\xi}(x)$. The $\nt_n$-modules $V_n^{\xi}$ are simple and form a natural chain whose union is $V_{\xi}$. Hence, $V_{\xi}$ is a simple $\nt$-module. We denote its annihilator in $U(\nt)$ by $J(f_{\xi})$.

\nota{Let $\Po(f_{\xi})$ be the set of all polarizations $\at$ of $\nt$ at $f_{\xi}$ such that $\at_n=\at\cap\nt_n$ is a polarization of $\nt_n$ at $f_n$ for large enough $n$. Define $V_{\xi,\at}$ and $V_{\xi,\at}^n$ by formula (\ref{formula:V_V_m_ifd}) in which $\pt$ and $\pt_n$ are replaced by $\at\in P(f_\xi)$ and $\at_n$ respectively. Then $V_{\xi,\at}=\ilm{V_{\xi,\at}^n}$ . The annihilator of $V_{\xi,\at}^n$ in $U(\nt_n)$ does not depend on $\at_n$, and \cite[Lemma 3.6]{IgnatyevPenkov1} and Proposition 3.6 below imply that the annihilator of $V_{\xi,\at}$ does not depend on~$\at$. This shows that $J(f_{\xi})$ can be defined via any polarization $\at\in\Po(f_{\xi})$.}

From now, we consider the cases $B_{\infty}$ and $D_{\infty}$. We may assume without loss of generality that the inclusions $\gt_n\subset\gt_{n+1}$ are defined by the natural embeddings $\Cp^{2n+1}\inj\Cp^{2n+3}$ (for $B_{\infty}$) and $\Cp^{2n}\inj\Cp^{2n+2}$ (for $D_{\infty}$) of the form $e_i\mapsto e_i$ (see Subsection~\ref{sst:center_fd}). Consequently, we may consider $\nt_n$ as a Lie subalgebra of $\nt_{n+1}$, and $U(\nt_n)$ as an associative subalgebra of $U(\nt_{n+1})$. Let $\Phi_n$ and $\Phi_{n+1}$ be the root systems of $\gt_n$ and $\gt_{n+1}$ respectively. Denote by $\Delta_n$, $\Delta_{n+1}$ the corresponding sets of simple roots, and put $\Bu_n=\Bu\cap\Phi_n$, $\Bu_{n+1}=\Bu\cap\Phi_{n+1}$. Given $i\in\{n,~n+1\}$, denote also
\begin{equation*}
\Bu_i'=\begin{cases}
(\Bu_i\setminus\Delta_i)\cup\{\epsi_{i-1}+\epsi_i\}&\text{if }\Phi_i=D_i\text{ with even }i,\\
\Bu_i\setminus\Delta_i&\text{otherwise}.
\end{cases}
\end{equation*}

By \cite{Dixmier4}, each $P_{\beta}$, $\beta\in\Bu$, acts on $V_{\xi}$ via some scalar $c_{\beta}$. It follows from \cite[6.6.9 (c)]{Dixmier3} that, given an arbitrary $\beta\in\Bu$,
$$c_{\beta}=\prod_{\beta'\in\Bu,~\col(\beta')\leq\col(\beta)}f(e_{\beta'})\neq0.$$
The following proposition plays the crucial role in the description of centrally generated primitive ideals of the algebra $U(\nt)$.

\propp{Each\label{prop:n_annihilating_n_plus_1} element of $U(\nt_n)$ annihilating $V_{\xi}^{n+1}$ belongs to the ideal of $U(\nt_{n+1})$ generated by all $P_{\beta}-c_{\beta}$ for $\beta\in\Bu_n'$.}{Put $\beta_1=\epsi_1+\epsi_2$, $e=e_{\epsi_1-\epsi_2}$, and, for $i=n,~n+1$, define $\kt_i$, $\theta_i$, $s_i$, $\Au_{s_i}$, $r_i$, $\wt\nt_e^i$, $\wt\nt_i$, $\wt\Phi_i$ for $\nt_i$ as $\kt$, $\theta$, $s$, $\Au_s$, $r$, $\wt\nt_e$, $\wt\nt$, $\wt\Phi$ respectively for $\nt$ in Subsection~\ref{sst:ideals_fd}. Note that $U(\wt\nt_e^i)=U(\wt\nt_i)\otimes\Cp[e]$. Recall that
$$r_i(e_{\epsi_1+\epsi_2}e_{\epsi_1-\epsi_2}+\ldots+e_{\epsi_1+\epsi_i}e_{\epsi_1-\epsi_i})=c_{\beta_1}(1\otimes e)\otimes1,~i=n,~n+1.$$
Since $\theta_{n+1}(x)=\theta_n(x)$ for $x\in\wt\nt_n$, the diagram
\begin{equation*}
\xymatrix{
U(\nt_n)\ar@{^{(}->}[d]\ar@{->>}[r]^-{r_n}&(U(\wt\nt_n)\otimes\Cp[e])\otimes\Au_{s_n}\ar[d]^{\vfi}\\
U(\nt_{n+1})\ar@{->>}[r]^-{r_{n+1}}&(U(\wt\nt_{n+1})\otimes\Cp[e])\otimes\Au_{s_{n+1}}
}
\end{equation*}
is commutative, where, by definition, $\vfi((x\otimes1)\otimes1)=(x\otimes1)\otimes1$ if $x\in\wt\nt_n$, $\vfi((1\otimes1)\otimes\overline{y})=(1\otimes1)\overline{y}$ if $y\in\kt_n$, and $$\vfi((1\otimes e)\otimes1)=(1\otimes e)\otimes1-c_{\beta_1}(1\otimes1)\otimes\overline{e}_{\epsi_1+\epsi_{n+1}}\overline{e}_{\epsi_1-\epsi_{n-1}}.$$

Now, denote by $J_i$ the annihilator of $V_{\xi}^i$ in $U(\nt_i)$, $i=n,~n+1$. Let $P_{\beta}^i$, $\beta\in\Bu_i$, be the canonical generators of $Z(\nt_i)$, $i=n,~n+1$. (Note that if $\row(\beta)<0$ then $P_{\beta}^i=P_{\beta}$.) According to Theorem~\ref{theo:ideal_fd_Bn_Dn} and \cite[6.6.9 (c)]{Dixmier3}, the ideal $J_i$ is generated by $P_{\beta}-c_{\beta}$, $\beta\in\Bu_i'$, and by $P_{\beta}^i$, $\beta\in\Bu_i\setminus\Bu_i'$. Now, put $\wt\Bu_i=\Bu_i\setminus\{\beta_1,~\epsi_1-\epsi_2\}$, $\wt\Bu_i'=\wt\Bu_i\cap\Bu_i'$, and denote by $\wt P_{\beta}^i$, $\beta\in\wt\Bu_i$, the canonical generators of $Z(\wt\nt_i)$, $i=n,~n+1$. It was shown in the proof of Proposition~\ref{prop:J_c_primitive_Bn_Dn} that there exist the unique nonzero scalars $a_{\beta}$, $\beta\in\wt\Bu_{n+1}$, such that $$r_{n+1}(P_{\beta}^{n+1})=a_{\beta}(\wt P_{\beta}^{n+1}\otimes1)\otimes1.$$ Since $\theta_{n+1}(x)=\theta_n(x)$ if $x\in\wt\nt_n$, one has $$r_n(P_{\beta}^n)=a_{\beta}(\wt P_{\beta}\otimes1)\otimes1,~\beta\in\wt\Bu_n',$$ where we denote $\wt P_{\beta}=\wt P_{\beta}^n=\wt P_{\beta}^{n+1}$. Furthermore, $r_i(J_i)=\wt J_e^i\otimes\Au_{s_i}$, where $\wt J_e^i$ is the ideal of $U(\wt\nt_i)\otimes\Cp[e]$ generated by $e$, $\wt P_{\beta}-\wt c_{\beta}$, $\beta\in\wt\Bu_i$, and by $\wt P_{\beta}^i$, $\beta\in\wt\Bu_i\setminus\wt\Bu_i'$ (here $\wt c_{\beta}=a_{\beta}^{-1}c_{\beta}$), $i=n,~n+1$.

Pick an element $x\in U(\nt_n)$ annihilating $V_{\xi}^{n+1}$. Let $\{\overline{a}_j\}$ be the basis of $\Au_{s_n}$ consisting of products of the powers of $\overline e_{\epsi_p\pm\epsi_q}$, $1\leq p\leq2$, $3\leq q\leq n$, taken in any fixed order. Then, clearly, there exist the unique $x_{j,k}\in U(\wt\nt_n)$ such that
$$r_n(x)=\sum_{j,k}(x_{j,k}\otimes e^k)\otimes\overline{a}_j=\sum_{j,k}((x_{j,k}\otimes1)\otimes\overline a_j)((1\otimes e)\otimes1)^k.$$
Applying $\vfi$ to this equality and denoting $\overline z=\overline{e}_{\epsi_1+\epsi_{n+1}}\overline{e}_{\epsi_1-\epsi_{n+1}}$, we have
\begin{equation*}
\begin{split}
r_{n+1}(x)&=\vfi(r_n(x))=\sum_{j,k}((x_{j,k}\otimes1)\otimes\overline a_j)(\vfi((1\otimes e)\otimes1))^k\\
&=\sum_{j,k}((x_{j,k}\otimes1)\otimes\overline{a}_j)((1\otimes e)\otimes1-c_{\beta_1}(1\otimes1)\otimes\overline z)^k=\sum_{j,k,l}(-1)^lc_{\beta_1}^l\binom{k}{l}(x_{j,k}\otimes e^{k-l})\otimes\overline z^l\overline a_j.
\end{split}
\end{equation*}
Here $\binom{k}{l}$ is the corresponding binomial coefficient, and $\binom{k}{l}=0$ if $k<l$.

Finally, denote by $\wt J_i$ the ideal of $U(\wt\nt_i)$ generated by $\wt P_{\beta}-\wt c_{\beta}$, $\beta\in\wt\Bu_i$, and by $\wt P_{\beta}^i$, $\beta\in\wt\Bu_i\setminus\wt\Bu_i'$, $i=n,~n+1$. Since $r_{n+1}(x)\in\wt J_e^{n+1}$ and $\overline z^l\overline a_j$ is a part of a basis of $\Au_{s_{n+1}}$, we see that $x_{j,k}\otimes e^{k-l}$ belongs to $\wt J_e^{n+1}$ for all $j,~k,~l$ with $l\leq k$. In particular, for $l=k$, we have that $x_{j,k}\otimes1\in\wt J_e^{n+1}$, or, in other words, $x_{j,k}\in\wt J_{n+1}$. But $x_{j,k}\in U(\wt\nt_n)$, thus we may use the induction on $n$ to show that in fact $x_{j,k}$ belongs to the ideal of $U(\wt\nt_n)$ generated by $\wt P_{\beta}-\wt c_{\beta}$ for $\beta\in\wt\Bu_n'$. (The base $n\leq2$ can be checked directly.) We conclude that $r_n(x)$ belongs to the ideal of $(U(\wt\nt_n)\otimes\Cp[e])\otimes\Au_{s_n}$ generated by $\wt P_{\beta}-\wt c_{\beta}$, $\beta\in\wt\Bu_n'$. Arguing as in the proof of Proposition~\ref{prop:J_c_primitive_Bn_Dn} again, we see that the preimage of this ideal is the ideal of $U(\nt_n)$ generated by $P_{\beta}-c_{\beta}$ for $\beta\in\Bu_n'$. The proof is complete.}

\corop{The\label{coro:J_f_xi_maximal_cent_gen} ideal $J(f_{\xi})$ is a maximal centrally generated ideal of $U(\nt)$.}{Pick an element $x\in J(f_{\xi})$. There exists $n$ such that $x\in U(\nt_n)$. On the other hand, $x$~annihilates $V_{\xi}^{n+1}\subset V_{\xi}$, hence, according to Proposition~\ref{prop:n_annihilating_n_plus_1}, $x$ belongs to the ideal of $U(\nt_n)$ generated by $P_{\beta}-c_{\beta}$ for certain $\beta\in\Bu$. Thus, $J(f_{\xi})$ is centrally generated.

Similarly, assume that $J$ is an ideal of $U(\nt)$ such that $J(f_{\xi})\subsetneq J$. Pick an element $x\in J$. There exists $n$ such that $x\in U(\nt_n)$, but $x$ annihilates $V_{\xi}^{n+1}$, hence, thanks to Proposition~\ref{prop:n_annihilating_n_plus_1}, $x\in J(f_{\xi})$, a~contradiction. Thus, $J(f_{\xi})$ is maximal.}

Our second main result is as follows (cf. \cite[Theorem 3.7]{IgnatyevPenkov1}).

\theop{Let $\Phi=B_{\infty}$ or $\Phi=D_{\infty}$\textup{,} and $\nt$ be as above. The following conditions on a primitive ideal $J\subset U(\nt)$ are equivalent:
\begin{equation*}
\begin{split}
&\text{\textup{i)} $J$ is centrally generated};\\[-2pt]
&\text{\textup{ii)} all scalars $c_{\beta}$, $\beta\in\Bu$, are nonzero};\\
&\text{\textup{iii)} $J=J(f_{\xi})$ for a Kostant form $f_{\xi}$}.\\
\end{split}
\end{equation*}
If these\label{theo:ideal_ifd} conditions are satisfied\textup{,} then the scalars $c_{\beta}$ reconstruct $\xi$ exactly as in Theorem~\textup{\ref{theo:ideal_fd_Bn_Dn}}.}{
$\mathrm{(ii)}\Longrightarrow\mathrm{(iii)}$. Define $\xi$ by formula (\ref{formula_xi_via_c_fd_Bn_Dn}). Then $J(f_{\xi})\subseteq J$ by Corollary~\ref{coro:J_f_xi_maximal_cent_gen}, but $J(f_{\xi})$ is maximal, so $J=J(f_{\xi})$.

$\mathrm{(iii)}\Longrightarrow\mathrm{(i)}$. Follows from Corollary~\ref{coro:J_f_xi_maximal_cent_gen}.

$\mathrm{(i)}\Longrightarrow\mathrm{(ii)}$. Assume, to the contrary, that some scalars $c_{\beta}$, $\beta\in\Bu$, equal zero. Suppose that $i_1$ is the minimal number such that $c_{\beta}=0$ for $\beta\in\Bu\cap\Ro_{i_1}$. We define inductively two (possibly, infinite) sequences $\{i_j\}$ and $\{k_j\}$ of positive integers by the following rule. If $i_j$ is already defined and there exists $k>i_j$ such that $c_{\beta}\neq0$ for $\beta\in\Bu\cap\Ro_k$, then set $k_j$ to be the minimal among all such $k$. Similarly, if $k_j$ is already defined and there exists $i>k_j$ such that $c_{\beta}=0$ for $\beta\in\Bu\cap\Ro_i$, then set $i_{j+1}$ to be the minimal among all such $i$.

To each $j$ for which both $i_j$ and $k_j$ exist we assign the set of roots
\begin{equation*}
\Gamma_j=\{\epsi_{i_j}+\epsi_{i_j+2},~\epsi_{i_j}-\epsi_{i_j+2}\}\cup\bigcup\limits_{i_j<s<k_j-1,~s\text{ even}}\{\epsi_s+\epsi_{s+3},~\epsi_s-\epsi_{s+3}\}\cup\{\epsi_{k_j-1}+\epsi_{k_j+1},~\epsi_{k_j-1}-\epsi_{k_j+1}\}.
\end{equation*}
If $i_j$ exists and $k_j$ does not exist, then put
\begin{equation*}
\Gamma_j=\{\epsi_{i_j}+\epsi_{i_j+2},~\epsi_{i_j}-\epsi_{i_j+2}\}\cup\bigcup\limits_{i_j<s,~s\text{ even}}\{\epsi_s+\epsi_{s+3},~\epsi_s-\epsi_{s+3}\}.
\end{equation*}
Now, we define the subset $X\subset\Phi^+$ as in the proof of Theorem~\ref{theo:ideal_fd_Bn_Dn}. Namely, denote the lengths of the sequences $\{i_j\}$, $\{k_j\}$ by $l_I$, $l_K$ respectively, and put
\begin{equation*}
X=\left(\Bu\cup\bigcup_{j=1}^{l_I}\Gamma_j\right)\setminus\left(\bigcup\{\Ro_i\mid c_{\beta}=0\text{ for }\beta\in\Bu\cap\Ro_i\}\cup\bigcup_{j=1}^{l_K}\Ro_{k_j}\right).
\end{equation*}

Let $$\mu_{\vfi}=\sum_{\alpha\in X}\vfi(\alpha)e_{\alpha}^*,$$ where $\vfi\colon X\to\Cp^{\times}$ is a map (note that this sum may be infinite). To each $\alpha=\epsi_i+\epsi_j\in\Phi^+$ we assign the subset $S(\alpha)\subset\Phi^+$ as follows:
\begin{equation*}
S(\alpha)=\bigcup_{k=i+1}^{j-1}\{\epsi_i-\epsi_k,~\epsi_i+\epsi_k,~\epsi_j+\epsi_k\}\cup\bigcup_{k=j+1}^n\{\epsi_i+\epsi_k\}\cup\Ro_j.
\end{equation*}
We then set $\Mo=\bigcup_{\beta\in X}\Mo_{\beta}$ where $$\Mo_{\beta}=\{\gamma\in S(\beta)\mid \gamma,\beta-\gamma\notin\bigcup\nolimits\Mo_{\alpha}\},$$ the latter union being taken over all $\alpha\in X$ such that $\row(\alpha)<\row(\beta)$. Note that if $\beta\in X$, $\alpha,\gamma\in\Phi^+$, $\alpha\notin\Mo$ and $\alpha+\gamma=\beta$, then $\gamma\in\Mo$. This implies that for all $x, y\in\at=\langle e_{\alpha},~\alpha\in\Phi^+\setminus\Mo\rangle_{\Cp}$ one has $\mu_{\vfi}([x,y])=0$. Moreover, it is easy to see that $\at$ is a subalgebra of $\nt$, hence we can consider the $\nt$-module $V_{\vfi}=U(\nt)\otimes_{U(\at)}W_{\vfi}$, where $W_{\vfi}$ is the one-dimensional representation of $\at$ given by $x\mapsto\mu_{\vfi}(x)$. Let $J_{\vfi}$ be the annihilator of $V_{\vfi}$ in $U(\nt)$ (we do not assert that $J_{\vfi}$ is a primitive ideal as we do not discuss the irreducibility of $V_{\vfi}$). One can check that $\vfi$ can be chosen so that $J\subseteq J_{\vfi}$, so we assume in the rest of the proof that this condition is satisfied.

Given $\gamma\in\Phi^+$, let $M$ be a finite subset of $\Zp_{>0}$ such that $\gamma\in\Phi_M^+$, $n=|M|$. Recall the definition of $j_M$ and $\phi_M$ from Subsection~\ref{sst:centers_ifd}. Let $\Delta_{\gamma}$, $P_{\gamma}$ be the respective images in $U(\nt_M)$ of $\Delta_{j_M^{-1}(\gamma)}$, $P_{j_M^{-1}(\gamma)}\in U(\nt_n)$ under the isomorphism $\phi_M$. Note that $\Delta_{\gamma}$ depends only on $\gamma$ and not on $M$. Let $\alpha=\epsi_{i_1}+\epsi_{i_1+2}$ be the unique root from $\Ro_{i_1}\cap X$. By Lemma~\ref{lemm:comm1_B_n_D_n}, $\Delta_{\alpha}$ commutes with all $e_{\gamma}$, except ${\gamma=\epsi_{i_1+1}-\epsi_{i_1+2}}$ and $\gamma=\epsi_{2j-1}+\epsi_{2j}$ for $j\geq1$, and $$[\Delta_{\alpha},e_{\epsi_{i_1+1}-\epsi_{i_1+2}}]=\pm\Delta_{\epsi_{i_1}+\epsi_{i_1+1}}.$$ Hence, $P_{\alpha}$ commutes with all~$e_{\gamma}$ except ${\gamma=\epsi_{i_1+1}-\epsi_{i_1+2}}$ and $\gamma=\epsi_{2j-1}+\epsi_{2j}$ for $j\geq1$,\break and $[P_{\alpha},e_{\epsi_{i_1+1}-\epsi_{i_1+2}}]=\pm P_{\epsi_{i_1}+\epsi_{i_1+1}}$. It follows from the finite-dimensional theory that $P_{\alpha}$ commutes with~$e_{\gamma}$ if $\row(\gamma)\neq i_1+1$ and $\col(\gamma)\neq\pm(i_1+1)$. Obviously, $P_{\alpha}$ commutes with $e_{\epsi_{i_1}-\epsi_{i_1+1}}$. Hence the image of~$P_{\alpha}$ belongs to the center of the image of $U(\nt)$ in the algebra $\End_{\Cp}V_{\vfi}$. By \cite{Dixmier4}, there exists $c\in\Cp$ such that $P_{\alpha}-c\in J$.

Further, we see that $$c=\pm\vfi(\alpha)\prod_{\beta\in\Bu,~\row(\beta)<i_1}\vfi(\beta),$$ because $J\subseteq J_{\vfi}$ and $\mu_{\vfi}(e_{\alpha})=\vfi(\alpha)$. Note also that, for any $\beta\in\Bu$, the element $P_{\beta}$ acts on $V_{\vfi}$ by the scalar $$c_{\beta}'=\prod_{\beta'\in\Bu,~\row(\beta')\leq\row(\beta)}\vfi(\beta').$$ Thus there exist at least two distinct maps $\vfi_1$, $\vfi_2$ from $X$ to $\Cp^{\times}$ such that $\vfi_1(\alpha)\neq\vfi_2(\alpha)$ and $c_{\beta}'=c_{\beta}$ for all $\beta\in\Bu$. This implies that both $J_{\vfi_1}$ and $J_{\vfi_2}$ contain $J$, which contradicts the uniqueness of $c$. The proof is complete.}\newpage

\nota{Denote by $\Au_{\infty}$ the Weyl algebra with countably many generators $p_i$, $q_i$ for $i\in\Zp_{>0}$ and relations
\begin{equation*}
[p_i,q_i]=1,~[p_i,q_j]=0\text{ for }i\neq j,~[p_i,p_j]=[q_i,q_j]=0\text{ for all }i,j.
\end{equation*}
It was proved in \cite[Corollary 3.8]{IgnatyevPenkov1} that, for $\Phi=A_{\infty}$ or $\Phi=C_{\infty}$ and $\Nu=\Zp_{>0}$, if $J$ is a centrally generated primitive ideal of $U(\nt)$ then the quotient algebra $U(\nt)/J$ is isomorphic to $\Au_{\infty}$. For $\Phi=B_{\infty}$ or $\Phi=D_{\infty}$ and $\Nu=\Zp_{>0}$, given a centrally generated primitive ideal $J$ of $U(\nt)$, denote by $\Au$\break the subalgebra of $U(\nt)/J$ generated by the images of $e_{\alpha}$, $\alpha\in\Phi^+\setminus\{\epsi_1-\epsi_2,~\epsi_3-\epsi_4,\ldots\}$,\break under the canonical projection $U(\nt)\to U(\nt)/J$. Using the elements $P_{\alpha}$ from the proof of Theorem~\ref{theo:ideal_ifd}, one can check that $\Au$ is isomorphic to~$\Au_{\infty}$. We do not know whether the entire quotient algebra $U(\nt)/J$ is isomorphic to~$\Au_{\infty}$.}

\medskip\textsc{Mikhail Ignatyev: Samara National Research University, Ak. Pavlova 1, 443011,\\\indent Samara, Russia}

\emph{E-mail address}: \texttt{mihail.ignatev@gmail.com}

\end{document}